\documentclass[11pt]{amsart}
\usepackage[left=1in,right=1in,top=1in,bottom=1in]{geometry}
\usepackage{graphicx}
\usepackage{amssymb}
\usepackage{epstopdf}
\newtheorem{theorem}{Theorem}

\newcommand{\csn}{\operatorname{CSN}}

\newcommand{\ofc}{\overline{\mathfrak{C}}}		
		
\newcommand{\pcsn}{\mathbb{P}\operatorname{CSN}}		
\newcommand{\bhv}{\operatorname{BHV}}		
\newcommand{\pbhv}{\mathbb{P}\operatorname{BHV}}		

\usepackage{graphicx}		
\usepackage{epstopdf}		
\usepackage{subcaption}		
\usepackage{float}		
\usepackage{wrapfig}		
\restylefloat{figure}	
\usepackage{amsthm}		
\usepackage{float}	
\usepackage{hyperref}	

\usepackage{tikz-cd}
		
\newcommand{\id}{\operatorname{id}}		
\newcommand{\Id}{\operatorname{Id}}
\newcommand{\mor}{\operatorname{Mor}}		
\newcommand{\ft}{\mathfrak{T}}		
\newcommand{\fc}{\mathfrak{C}}		
\newcommand{\complex}{\mathbb{C}}		
\newcommand{\real}{\mathbb{R}}

\newcommand{\gl}{\mathrm{GL}}
\newcommand{\pgl}{\mathrm{PGL}}
\newcommand{\sgl}{\mathrm{SL}}
\newcommand{\psl}{\mathrm{PSL}}
\newcommand{\cp}{\mathbf{C}\mathbb{P}}
\newcommand{\rp}{\mathbf{R}\mathbb{P}}
\newcommand{\dih}{\mathrm{Dih}}

\usepackage[colorinlistoftodos]{todonotes}
\usepackage{algorithm}
\usepackage{tikz}
\usetikzlibrary{matrix}
\usepackage{xcolor}

\usepackage[noend]{algpseudocode}

\newtheorem{definition}[theorem]{Definition}
\newtheorem{proposition}[theorem]{Proposition}
\newtheorem{corollary}[theorem]{Corollary}
\newtheorem{lemma}[theorem]{Lemma}

\title{Comparison Theorems of Phylogenetic Spaces and the Moduli Spaces of Curves}
\date{\today}

\author{Y. Wu}
\author{S.-T. Yau}

\begin{document}
\begin{abstract}
Rapid developments in genetics and biology have led to phylogenetic methods becoming an important direction in the study of cancer and viral evolution. Although our understanding of gene biology and biochemistry has increased and is increasing at a remarkable rate, the theoretical models of genetic evolution still use the phylogenetic tree model that was introduced by Darwin in 1859 and the generalization to phylogenetic networks introduced by Grant in 1971. Darwin's model uses phylogenetic trees to capture the evolutionary relationships of reproducing individuals \cite{darwin1859origin}; Grant's generalization to phylogenetic networks is meant to account for the phenomena of horizontal gene transfer \cite{grant1971}. Therefore, it is important to provide an accurate mathematical description of these models and to understand their connection with other fields of mathematics. In this article, we focus on the graph theoretical aspects of phylogenetic trees and networks and their connection to stable curves. We introduce the building blocks of evolutionary moduli spaces, the dual intersection complex of the moduli spaces of stable curves, and the categorical relationship between the phylogenetic spaces and stable curves in $\overline{\mathfrak{M}}_{0,n}(\mathbb{C})$ and $\overline{\mathfrak{M}}_{0,n}(\mathbb{R})$. We also show that the space of network topologies maps injectively into the boundary of $\overline{\mathfrak{M}}_{g,n}(\mathbb{C})$.
\end{abstract}
\maketitle
\tableofcontents

\section{Introduction}\label{sec:intro}

The evolutionary process is usually described by a rooted phylogenetic tree to characterize a set of species, the spread of a pathogen, or a tumor mutation \cite{huson2010phylogenetic, zairis2014moduli}. Phylogenetic trees are well suited to represent evolutionary histories, where the main events are speciations represented as the internal nodes and descent with modifications denoted along the edges of the tree. The mathematical model of the space of phylogenetic trees was introduced by Billera-Holmes-Vogtmann using the metric geometry of tree spaces \cite{billera2001geometry}.

A striking fact revealed by recent phylogenetic analysis is that many eukaryotic nuclear genes are of bacterial origin (i.e., prokaryotic ancestry). It is now believed that this is due to the replacement of nuclear genes of early eukaryotes with those of the prokaryotes that they engulfed for food \cite{doolittle1998you}. Such gene transfer was called horizontal gene transfer (HGT), which more generally refers to the transfer of DNA between organisms to permit the acquisition of novel traits \cite{hotopp2011horizontal,jain2003horizontal}. Reticulate evolution may enhance adaptive radiation by introducing genetic variation, enabling the population to take advantage of new ecological opportunities through biodiversity \cite{lamichhaney2015evolution,pease2016phylogenomics}. Therefore, reticulate evolution across the genes has been found to be extensive \cite{Edelman594,kozak2018genome}. Phylogenetic networks have been introduced to model this remarkable phenomenon, along with many other reticulate events, such as hybridization, recombination, gene conversion, or gene duplication \cite{huson2010phylogenetic}. The space of phylogenetic networks was interpreted by Devadoss and Petit  as a natural extension of the space of phylogenetic trees \cite{devadoss2017space}.

The space of phylogenetic trees and the space of phylogenetic networks are metric spaces with simplicial structures, and we call them phylogenetic moduli spaces.
The geometry of the spaces of phylogenetic trees and networks was constructed to accurately model what is known about the transmission of genetic mutations from one generation to the next (from parent to child) and across generations (by the mechanism of horizontal transfer).  This is to say that the definitions and structure of the spaces of phylogenetic trees and networks are mandated by the biology of genes and gene transfer.  Their definitions are biologically canonical.  Thus it is a surprising fact that various versions of these spaces are instances of canonical mathematical spaces such as simplicial complexes \cite{wu2019comparison},  cubical complexes \cite{billera2001geometry}, tropical Grassmannians \cite{speyer2004tropical}, and algebraic fans spanned over simplicial complexes formed by root systems of type $D$ \cite{wu2019comparison}. There is also a metric on the spaces of phylogenetic trees and phylogenetic networks \cite{ardila2012geodesics}, and Owen et al.~proved that geodesics in spaces of phylogenetic trees can be computed in polynomial time \cite{owen2011computing,owen2011fast}. From an algebraic geometry approach, Devadoss and Morava formulated a smooth blowup to spaces of phylogenetic trees from the compactified moduli stacks of stable genus $0$ curves with $n$ marked points $\mathfrak{M}_{0,n}(\mathbb{R})$ \cite{devadoss2010diagonalizing}; from an algebraic topology approach, Baez and Otter constructed an operad whose operations are the edge-labeled phylogenetic trees \cite{baez2017operads}. 
 
In this paper, we discuss the connection between the phylogenetic moduli spaces and the moduli spaces of stable curves. These two subjects are seemingly distant; the moduli spaces of curves are fundamental objects of mathematics, whereas the phylogenetic spaces are canonical biological objects.  They have no prior knowledge of each other.  But, as is explained in this paper, these two spaces, defined in distant fields of science, are closely related by virtue of having equivalent categorical structures.  
Whether this close relationship can be exploited for novel biological applications remains to be seen.  In any event, this sort of close relationship was in no way expected.

\begin{figure}[!htb]
 \centering
 \includegraphics[width=.95\textwidth]{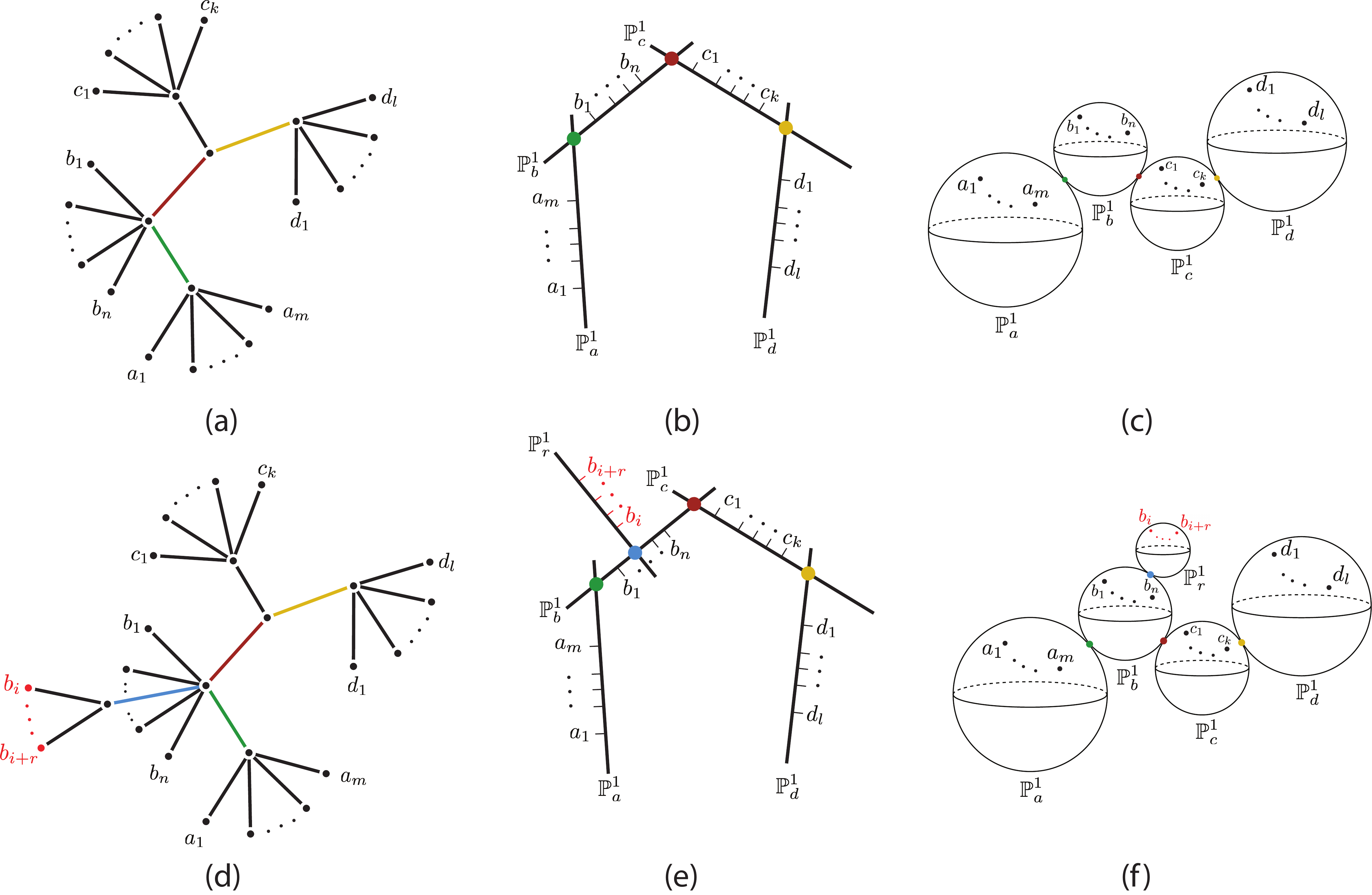}
\caption{Phylogenetic Trees (left) and the corresponding boundary stratum in the space of stable curves $\overline{\mathfrak{M}}_{0,n}(\mathbb{C})$ (middle and right).}
\label{fig:stable}\end{figure}

\vskip 3mm

\noindent\bf Acknowledgements. \rm The authors would like to thank Cliff Taubes
for a great deal of inspiration from his work and many fruitful discussions regarding the content of this paper.

\section{Algebraic Morphisms Between Phylogenetic Trees and Stable Curves in $\overline{\mathfrak{M}}_{0,n}(\mathbb{C})$}
In this section, we introduce the building blocks of evolutionary moduli spaces and the dual intersection complex of the moduli spaces of stable curves. Then we introduce comparison theorems between the categories of these spaces. 
Figure \ref{fig:stable} (a)-(c) gives an example of the correspondence between a generic phylogenetic tree and a boundary stratum in the space of stable curves $\overline{\mathfrak{M}}_{0,n}(\mathbb{C})$. 
Figure \ref{fig:stable} (d)-(f) illustrates extending the $b_i, \ldots, b_{i+r}$ nodes from the tree and the corresponding changes in the respective boundary stratum. More detailed correspondence between genus 0 stable curves and acyclic graphs is given in Section \ref{subsec:trees} for the genus 0 case after necessary definitions are introduced, and the generic correspondence between high genus stable curves and graphs introduced by Arbarello, Cornalba 􏱥and Griffiths \cite{arbarello2011geometry} is presented in Section \ref{subsec:networks}. Note that we need at least three leaves plus a root to define a phylogenetic tree and four leaves to define a phylogenetic network. Also, we need at least four marked points to have nonempty boundary for $\overline{\mathfrak{M}}_{0,n}$. The leaves in phylogenetic trees and phylogenetic networks correspond to marked points on algebraic curves, as we discuss in detail in this paper. Therefore, we restrict our discussion to at least three leaves plus a root for phylogenetic trees and at least four leaves for phylogenetic networks, corresponding to $\overline{\mathfrak{M}}_{0,n}$ with $n \geq 4$ throughout this paper.

\subsection{The Space of Phylogenetic Trees}\label{subsec:trees}
We use $\mathcal{X} = \{x_1, \ldots, x_n\}$ to denote a set of taxa, in which each taxon $x_i$ represents some species, group or individual organism whose evolutionary history is of interest to us \cite{huson2010phylogenetic}. For example, we can have $\mathcal{X}$ denote a set of mammals, with $x_1$ representing gorillas, $x_2$ representing seals, and so on. For our purposes, the set $\mathcal{X}$ of taxa is a set labeled by integers from 1 up to some maximal integer (which we call $n$). A phylogenetic tree with $n$ leaves is the equivalence class of a weighted, connected graph with no cycles, having $n$ distinguished labeled vertices up to rotations and a root. The labeled vertices and the root are of degree 1 and are called leaves. All the other vertices are of degree $\geq$ 3.
 A formal definition is presented as the following:

\begin{definition}\rm\label{def:tree}
A phylogenetic tree $T^{\mathcal{X}}$ is an acyclic graph, with internal edges being edges with two vertices and $n$ external edges being edges with just a single vertex.  The external edges are in 1-1 correspondence with the set $\mathcal{X}$ and are called the \it leaves of the graph. \rm Meanwhile, each internal edge has an associated \it length, \rm which is a positive real number.  The set of lengths of the internal edges is denoted by $\mathcal{W}^\mathcal{X}$.  The graph with the internal edge weight assignment is constrained by its set of internal edges that can be decomposed as a union of disjoint subsets of edges, called \it splits, \rm such that each edge separates the set $\mathcal{X}$ as the union of two disjoint subsets with both subsets having at least two leaves. For present purposes, a split is the same as an internal edge; we call it a split because if we remove it, the graph becomes a disjoint union of two trees.
The set of splits is denoted by $\mathcal{S}^\mathcal{X}$; it is an ordered set with components $\left\{S_1^{A_1| B_1}, \ldots, S_k^{A_k | B_k}\right\}$ with any given pair $(A_i, B_i)$ being the corresponding disjoint subset decomposition of $\mathcal{X}$ with each split a priori determined by the decomposition of $\mathcal{X}$.  The set $\mathcal{W}^\mathcal{X}$ of weights of the internal edges is sometimes written as $\left\{w_1^{A_1| B_1}, \ldots, w_k^{A_k | B_k}\right\}$ with each constituent being the weight assigned to the edges in the like-labeled split. 
 When two distinct phylogenetic trees have the same set $\mathcal{X}$, we denote them $T^\mathcal{X}_x$ and $T^\mathcal{X}_y$. 
\end{definition}

A phylogenetic tree contains the information of the set of bipartitions $\mathcal{S}^\mathcal{X}_x$ and their length $\mathcal{W}^\mathcal{X}_x$; therefore, we write it as 
$T^\mathcal{X}_x(\mathcal{S}^\mathcal{X}_x, \mathcal{W}^\mathcal{X}_x)$. If phylogenetic trees $T^\mathcal{X}_x$ and $T^\mathcal{X}_y$ induce the same set of bipartitions, we say that they have the same \it tree topology, \rm denoted as $\overline{T}^\mathcal{X}_x\sim \overline{T}^\mathcal{X}_y$. We use $\mathfrak{W}^z_n$ to denote the moduli space of the set of edge lengths of the tree topology $\overline{T}^\mathcal{X}_z$, with interior points homeomorphic to $\real^k_{>0}$, with $k$ being the number of internal edges of the tree. In other words, the tree topology forgets the length of an internal edge but delineates whether each internal edge is equal to zero. If $\mathcal{S}^\mathcal{X}_x \setminus \mathcal{S}^\mathcal{X}_y \neq \emptyset$, then we say the trees have different topology.

\begin{figure}[!htb]
 \centering
 \includegraphics[width=.7\textwidth]{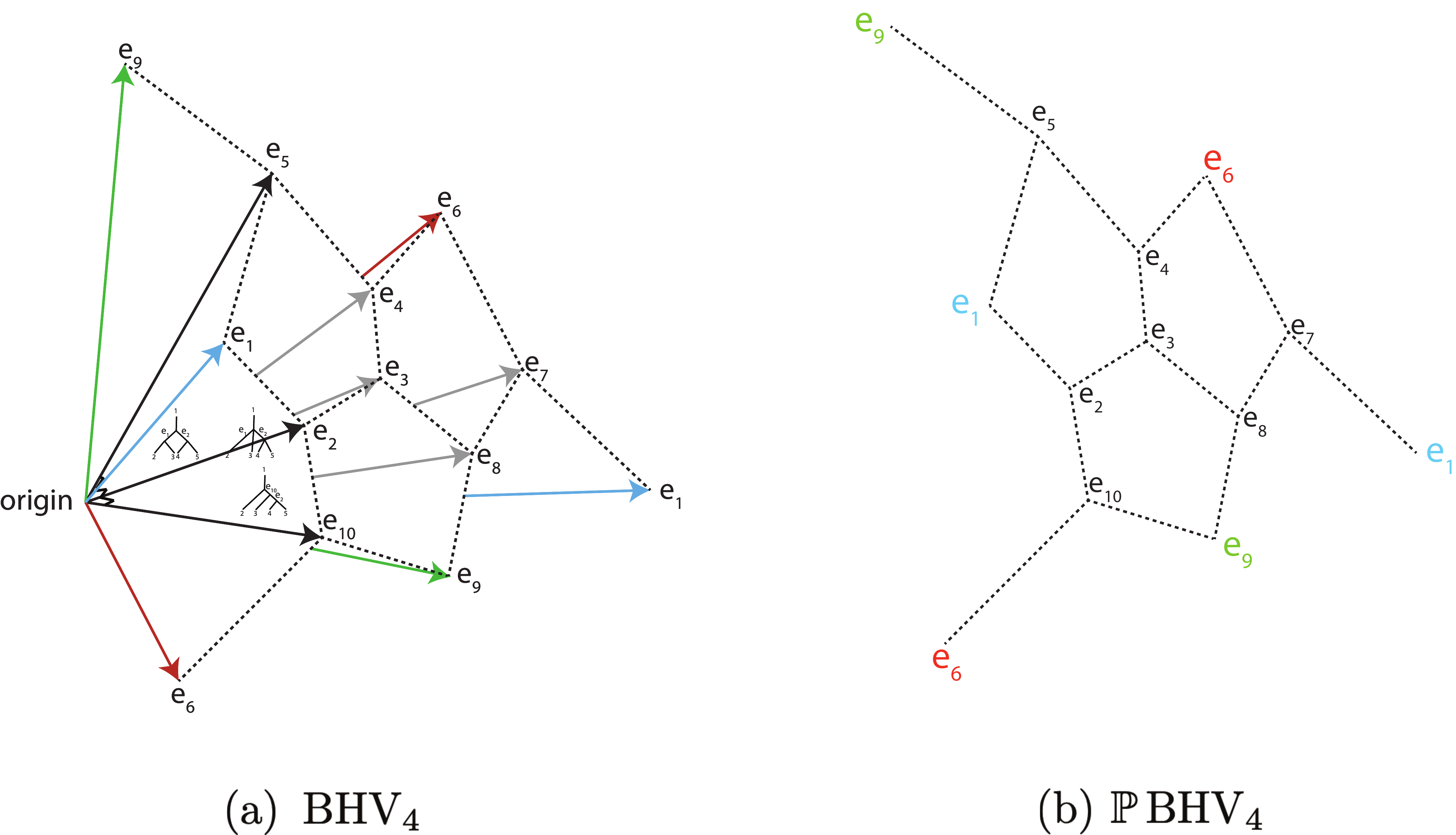}
 \caption{The space of phylogenetic trees and the projective space of phylogenetic trees.}
 \label{fig:bhvpbhv}
\end{figure}

\begin{definition}\label{def:BHV}
The \it Billera-Holmes-Vogtmann space of phylogenetic trees \rm is the space of isometry classes of rooted phylogenetic trees with $n$-labeled leaves, where the nonzero weights are on the internal branches, denoted as $\bhv_n$. The space $\bhv_n$ is constructed by gluing together $(2n - 3)!!$ positive orthants; each orthant corresponds to a particular tree topology, with the coordinates specifying the lengths of the edges. A point $T^\mathcal{X}_x$ in the interior of an orthant in $\bhv_n$ represents a binary tree, which is a tree in which each node has at most two children. If any of the coordinates is 0, the tree is called a \it degenerated tree, \rm which is obtained by collapsing edges corresponding to the 0 coordinates from a binary tree. We glue orthants together, such that a tree is on the boundary between two orthants when it can be obtained by collapsing edges from either tree topology. The $\bhv_n$ space is a cubical complex built by tiling each orthant with unit cubes of dimension $n-2$ \cite{zairis2016genomic}. 
\end{definition}

\begin{definition}\label{def:PBHV}
The \it projective Billera-Holmes-Vogtmann space of phylogenetic trees
$\pbhv_n$ \rm is the subspace of $\bhv_n$ consisting of points in each orthant, for which the sum of its $n-2$ internal edges is 1.
$\pbhv_n$ inherits a simplicial structure from $\bhv_n$, where the $k$-simplices of $\pbhv_n$ are points projected from points in $\bhv_n$ with exactly $k+1$ nonzero edges, and the intersection of faces is determined by degenerated trees that share the same internal edges. 
\end{definition}

$\bhv_4$ and $\pbhv_4$ spaces are illustrated in Figure \ref{fig:bhvpbhv}, where edges in \ref{fig:bhvpbhv}(a) and vertices in \ref{fig:bhvpbhv}(b) are identified by color, i.e., the two copies of $e_1$, $e_6$, and $e_9$ are identified. The structure of the internal branches is captured by the $\bhv_{n-1}$ construction. Zairis, Khiabanian, Blumberg, and Rabadan constructed the \it evolutionary moduli space \rm \cite{zairis2014moduli}, which allows potentially nonzero weights for the $n$ external leaves by crossing with an $n$-dimensional orthant:
$$\Sigma_n = \bhv_{n-1} \times \mathbb{R}_{\geq 0}^n.$$

In the following, we introduce relevant notions in order to show the connection between phylogenetic trees and stable curves in $\mathfrak{M}_{0,n}(\mathbb{C})$.

\begin{definition}\rm\label{def:nodal}
A nodal curve is a complete algebraic curve such that every one of its points is either smooth or is locally complex-analytically isomorphic to a neighborhood of the origin in the locus with equation $xy = 0$ in $\mathbb{C}^2$ \cite{arbarello2011geometry}.\end{definition}

\begin{definition}\rm\label{def:stable}
An $n$-pointed \it nodal curve \rm consists of the datum $(C; p_1, \ldots, p_n)$ of a nodal curve $C$ together with $n$ distinct smooth points of $C$. 
Let $C$ be a nodal curve, and let $D$ be a finite set of smooth points of $C$. Let $(C;D)$ be a connected nodal curve with $n$ marked points. 
$(C;D)$ is said to be \it stable \rm if it has a finite automorphism group \cite{arbarello2011geometry}.
\end{definition}

\begin{definition}\rm
The \it moduli space of an $n$-pointed stable curve of given genus $g$ \rm is denoted by $\overline{\frak{M}}_{g,n}(\mathbb{C})$. As a set, it is the set of isomorphism classes of $n$-pointed genus $g$ stable curves. Denote by $\frak{M}_{g,n}(\mathbb{C})$ the set of isomorphism classes of \it smooth \rm $n$-pointed genus $g$ stable curves \cite{arbarello2011geometry}. 
\end{definition}

\begin{definition}\label{defn:compact}\rm
$\overline{\mathfrak{M}}_{g,n}(\mathbb{C})$ \it compactifies \rm $\mathfrak{M}_{g,n}(\mathbb{C})$ without ever allowing the points to come together. When points on a smooth curve approach each other, the curve sprouts off one or more components, each isomorphic to the projective line \cite{fulton1996notes}.
\end{definition}

This compactification is called the Deligne-Mumford-Knudsen compactification, constructed by Deligne and Mumford \cite{deligne1969irreducibility} for $n = 0$ and by Knudsen \cite{deligne1972irreducibility} in general. 

\begin{definition}\rm
A nonempty subset $Y$ of a topological space $X$ with Zariski topology is \it irreducible \rm if it cannot be expressed as the union $Y = Y_1 \cup Y_2$ of two proper subsets, each of which is closed in $Y$. The empty set is not considered to be irreducible \cite{hartshorne2013algebraic}.
\end{definition}

\begin{definition}\rm\label{def:divisor}
A \it divisor \rm is an element of the free abelian group generated by the subvarieties of codimension one \cite{hartshorne2013algebraic}.
\end{definition}

\begin{definition}\rm\label{def:boundary}
The \it boundary \rm $\Delta = \overline{\mathfrak{M}}_{g,n}(\mathbb{C}) - \mathfrak{M}_{g,n}(\mathbb{C})$ is a divisor, with each component the closure of a locus of curves with 1 node \cite{harris2006moduli}.
\end{definition}

\begin{definition}\rm\label{fulton1996notes}
The \it boundary divisor \rm of $\overline{\mathfrak{M}}_{0,n}$ corresponding to the
marking partition $A \cup B = [n]$ is naturally isomorphic (by gluing) to the product
$$\overline{\mathfrak{M}}_{0,A\cup\{\cdot\}} \times \overline{\mathfrak{M}}_{0,B\cup\{\cdot\}}.$$
\end{definition}

A stable curve with one node is either irreducible, or is the union of smooth curves of genera $i$ and $g - i$ meeting at one point \cite{harris2006moduli}. We can count the number of complex parameters in these two components, and the dimension of this stratum is their sum. The first component with $k$ marked points has $k+1-3$ complex parameters and the second component with $n-k$ marked points has $n-k+1-3$ complex parameters. Therefore, they add up to be dimension $n-4$. Since $\overline{\mathfrak{M}}_{g,n}(\mathbb{C})$ has complex dimension $3g+n-3$ 
 \cite{fulton1996notes}, $\overline{\mathfrak{M}}_{0,n}(\mathbb{C})$ has complex dimension $n-3$. Therefore, these strata have complex codimension 1, and they are boundary divisors \cite{fulton1996notes}.

\begin{figure}[!htb]\centering
 \includegraphics[width=.75\textwidth]{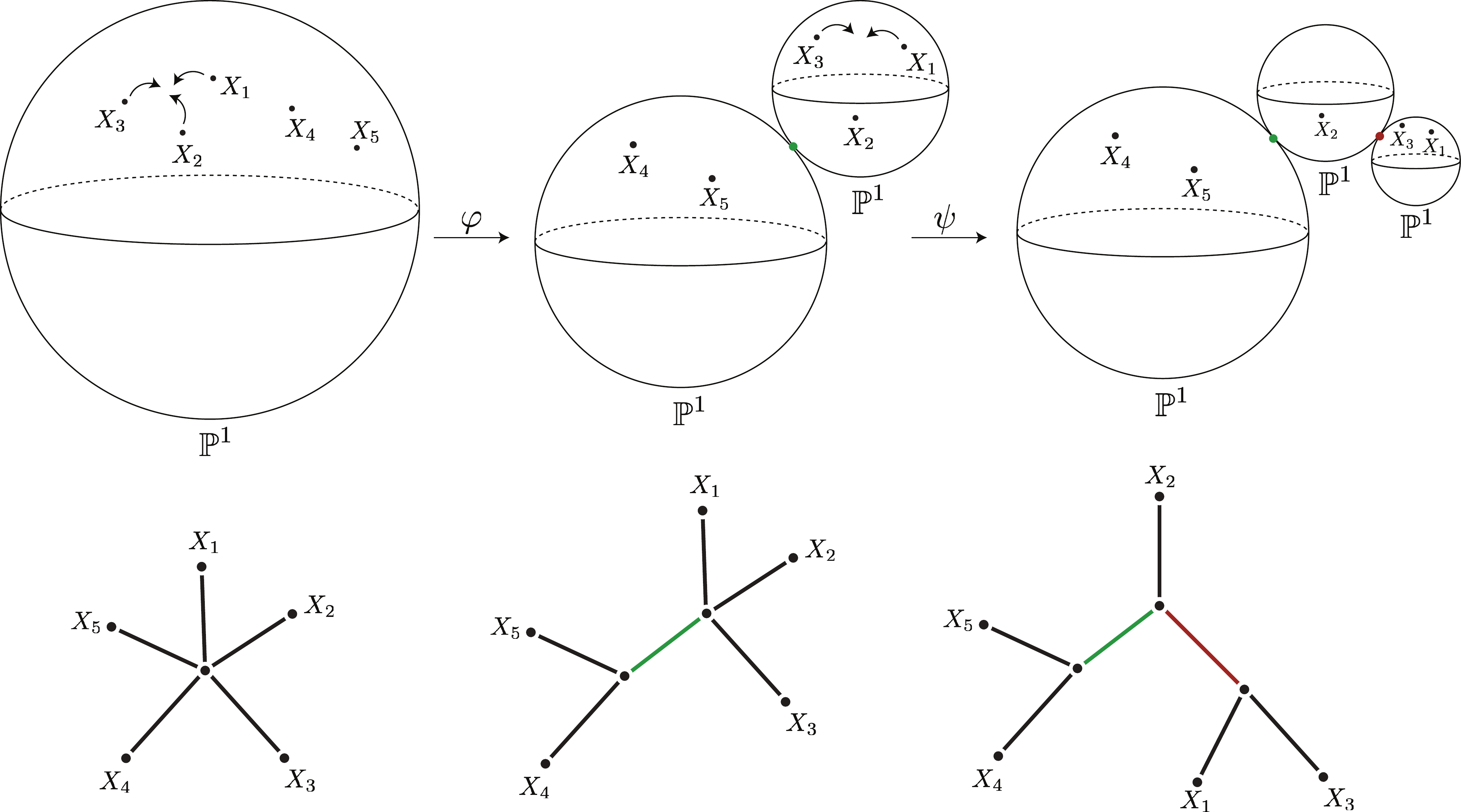}
\caption{Collision of three points in $\mathfrak{M}_{0,5}(\mathbb{C})$. \it Left: \rm an interior point of $\mathfrak{M}_{0,5}(\mathbb{C})$ (top) and its corresponding tree (bottom); \it middle: \rm a boundary divisor of $\overline{\mathfrak{M}}_{0,5}(\mathbb{C})$ with one collision (top) and its corresponding tree (bottom); 
\it right\rm: a boundary on $\overline{\mathfrak{M}}_{0,5}(\mathbb{C})$ with two collisions (top) and its corresponding tree (bottom).}\label{fig:collide}
\end{figure}

Figure \ref{fig:collide} depicts a collision of three points in $\mathfrak{M}_{0,5}(\mathbb{C})$. The left picture is an interior point of $\mathfrak{M}_{0,5}(\mathbb{C})$. Because there is a M\"obius group acting on these five marked points, there are two complex parameters for this interior point of $\mathfrak{M}_{0,5}(\mathbb{C})$. Colliding three marked points on an interior point of $\mathfrak{M}_{0,5}(\mathbb{C})$ sprouts another $\mathbb{P}^1$ with these 3 points, as displayed in the middle picture of Figure \ref{fig:collide}, ending one complex parameter on the $\mathbb{P}^1$ with these three marked points and having no complex parameter on the component with two marked points. This shows that colliding reduces the degree of freedom by one complex dimension. We can continue to collide two points on the $\mathbb{P}^1$ with three marked points in the middle, yielding the picture on the right in Figure \ref{fig:collide}.

The bottom row of Figure \ref{fig:collide} depicts corresponding trees. For each intersection, we denote the set of marked points on one side as $A_i = \{X_{i_1}, \ldots, X_{i_k}\}$ and on the other side as $B_i = \{X_{i_{k+1}},\ldots, X_{i_n}\}$. This intersection corresponds to an internal edge of a tree representing the edge $S_n^{A_i|B_i}$. For example, the stable curve in the middle has two components with the set of marked points $A_1 = \{X_1, X_2, X_3\}$ and $B_1 = \{X_4, X_5\}$, respectively. This stable curve corresponds to a tree with an internal edge $S_5^{A_1|B_1}$. After the second collision, the stable curve to the right has three components; for the first intersection, the set of marked points on one side is $A_1$ and on the other side is $B_1$. For the second intersection, the set of marked points on one side is $A_2 = \{X_2, X_4, X_5\}$ and $B_2 = \{X_1, X_3\}$; this stable curve corresponds to a tree with two internal edges $S_5^{A_1|B_1}$ and $S_5^{A_2|B_2}$.

\subsection{Cellular Decomposition of $\overline{\mathfrak{M}}_{0,n}(\mathbb{R})$} 
In this section, we introduce the cellular decomposition of $\overline{\mathfrak{M}}_{0,n}(\mathbb{R})$ based on Keel's construction of the moduli space \cite{keel1992intersection} and discussions in Harris and Morrison \cite{harris2006moduli}.
The moduli space of algebraic curves with $n$ marked points $\overline{\mathfrak{M}}_{0,5}(\mathbb{C})$ can be regarded as the quotient of the space of $n$-marked points on $\cp^1$ by the natural action of the group $\mathrm{SL}(2; \complex)$. Thus, the quotient is an equivalent class up to M\"obius transformations. The group $\gl(2; \complex)$ acts on $\complex \cup \{\infty\}$, which is $S^2$ topologically:
\begin{align*}
\gl(2; \complex) \ &\text{\rotatebox[origin=c]{-90}{$\circlearrowright$}}\ \complex \cup \{\infty\}\\
\begin{pmatrix}
a&b\\
c&d\end{pmatrix}&\cdot x = {ax + b\over cx + d}.
\end{align*}
Since $x$ is sent to this same point for
$rM=
\begin{pmatrix}
ra&rb\\
rc&rd\end{pmatrix}$ for all $r \in \mathbb{C} \setminus \{0\}$, the action factors through $\pgl(2; \complex)$. Without loss of generality, we assume that the automorphism group fixes the first three points $X_1$, $X_2$, and $X_3$ at $0, 1$, and $\infty$, respectively. Then we have two points, $X_4$ and $X_5$, moving on $\cp^1$, and we have the constraint that all points are distinct, i.e., (1) $X_4$ and $X_5$ cannot be the same as $X_1$, $X_2$, and $X_3$; (2) $X_4$ and $X_5$ cannot be the same as each other. From the first constraint, we have the moduli space of $X_4$ and $X_5$ being a product of two copies of $\cp^1 \setminus \{0,1,\infty\}$. For the second constraint, we need to take out the diagonal, which is the subspace with $X_4 = X_5$. Therefore, we have
$$\overline{\mathfrak{M}}_{0,5}(\mathbb{C}) = \underbrace{(\cp^1 \setminus \{0,1,\infty\}) \times \cdots \times (\cp^1 \setminus \{0,1,\infty\})}_\text{$n-3$ factors} \setminus \Delta,$$
where $\Delta$ is the thick diagonal \cite{keel1992intersection}.

The moduli space of algebraic curves over real numbers $\overline{\mathfrak{M}}_{0,n}(\mathbb{R})$ is
$$\overline{\mathfrak{M}}_{0,5}(\mathbb{R}) = \underbrace{(\rp^1 \setminus \{0,1,\infty\}) \times \cdots \times (\rp^1 \setminus \{0,1,\infty\})}_\text{$n-3$ factors} \setminus \Delta$$
with $\pgl(2; \real)$ acting on it.
Since $\rp^1$ is a circle, the marked points have a cyclic order. Therefore, we fix the first three points at $0, 1$, and $\infty$ by the automorphism group $\pgl(2; \real)$.
Also, the coordinates of $X_k\in\real$ are related to the $[0, 2\pi)$ coordinate $t_k$ via the rule whereby
$$X_k = {i (1 -e^{it_k})\over 1+ e^{it_k}}= \tan {t_k \over 2}.$$
Thus, the following $n-3$ points move between $(\pi, 2\pi)$ on $\rp^1$.

The configuration of 5 distinct points on the circle $\overline{\mathfrak{M}}_{0,5}(\mathbb{R})$ forms a union of open cells. Because $\mathfrak{M}_{0,n}(\mathbb{R})$ is defined without allowing marked points to come together, the boundary faces in
$\partial \overline{\mathfrak{M}}_{0,n}(\mathbb{R}) =\overline{\mathfrak{M}}_{0,n}(\mathbb{R}) - \mathfrak{M}_{0,n}(\mathbb{R})$
 arise from the points colliding.
We consider the cell with the ordering $(X_1, X_2, X_3, X_4, X_5)$, such that $X_1$, $X_2$, and $X_3$ are located at $t = 0, {\pi\over 2}$, and $\pi$. These three points represent $0,1$, and $\infty$ compactified on $\rp^1$; Meanwhile, $X_4,X_5$ range between $(\pi, 2\pi)$. If we place point $X_1$ at $t_1 = 0$, then each of the other points has an interval in which to move. Moving from 0 towards $2\pi$, we can use the residual $\psl(2; \real)$ action to make the next point after $X_1$ fixed at $t = 1$ and the next point after that fixed at $t = \pi$ (take that to be $\infty$).  The last two points move between $(\pi, 2\pi)$.  The parameter space for this is a 2-dimensional simplex, with one parameter $X_4$ moving between $(\pi, 2\pi)$ and $X_5$ moving in $(X_4, 2\pi)$.
Since we set $X_1$ at $t = 0$, the number of components is the number of ways to order the 4 remaining points, i.e., $4! = 24$ ways under $\psl(2; \real)$ and $24/2=12$ ways under $\pgl(2; \real)$. These configurations form 12 associahedra tiling $\overline{\mathfrak{M}}_{0,5}(\mathbb{R})$, and each associahedra is called a ``cell'' \cite{devadoss2017space}. $\overline{\mathfrak{M}}_{0,6}(\mathbb{R}) $ is tiled by 60 associahedra. In general, possible arrangements of $n$ points on the circle $\rp^1$ are:
$${n! \over  |\dih_n|} = 60,$$
where $|\dih_n|=2n$ represents the dihedral group, and the quotient identifies the rotation and reflection of an $n$-gon.

 Any matrix $M \in \pgl(2; \real)$ can be written as $M = AT$ with $A \in \sgl(2; \real)$, $T$ acting on $\real \cup \{\infty\}$ sending $x$ to $-x$ when we view the point as a real number, and $T$ changing $x$ to $2\pi - t$ when we view the point on the compactification of $\mathbb{R}$ as $S^1$:
 \begin{align*}
T: \rp^1 & \longrightarrow \rp^1\\
x & \mapsto \begin{pmatrix}-1&0\\0&1\end{pmatrix}x.
 \end{align*}
 Therefore, when the matrix $T$ is viewed as a transformation on $[0, 2\pi)$, it reverses the cyclic order of a set of ordered points on $[0, 2\pi)$, i.e., it sends a cyclic order $X_1$, $X_2$, $X_3$, $X_4$, $X_5$ to the cyclic order $X_5$, $X_4$, $X_3$, $X_2$, $X_1$, since if $t_1 < t_2 < t_3 < t_4 < t_5$, then $2\pi - t_5 < 2\pi - t_4 < 2\pi - t_3 < 2\pi - t_2 < 2\pi - t_1$.  Note that $T$ switches the cyclic order, rather than mapping $X_1,X_2,X_3,X_4$, and $X_5$ to $X_5,X_4,X_3,X_2$, and $X_1$ respectively. Then we can use $T$ to flip $(X_1,X_2,X_3,X_4,X_5)$ into $(X_5,X_4,X_3,X_2,X_1)$, and then use $A\in\sgl(2; \real)$ to write $(X_5,X_4,X_3,X_2,X_1)$ as $(X_1, X_5,X_4,X_3,X_2)$ with $X_1$ fixed at $t_1 = 0$, $X_5$ fixed at $t_5 = {\pi \over 2}$, and $X_4$ fixed at $t_4 = \pi$. This leaves $X_3$ and $X_2$ free to move in $(\pi, 2\pi)$ with $X_2$ having a $t$-coordinate greater than $X_3$.

$X_5$ colliding with $X_1$ results in a boundary stratum consisting of one $\rp^1$ with three special points $(X_1, x, X_5)$ at $\left(0, {\pi\over2}, \pi\right)$ and a second $\rp^1$ (the original one) with points $(x, X_2, X_3, X_4)$. When we look at boundary strata, separate elements of $\pgl(2;\real)$ act on each component $\rp^1$. This face is 1-dimensional because only $X_4$ on the second component is allowed to move in $(\pi, 2\pi)$, $X_1, x, X_5$ on the first $\rp^1$ are fixed at $\left(0, {\pi\over2}, \pi\right)$ by the group action $\pgl(2;\real)$, and $x, X_2, X_3$ are fixed at $\left(0, {\pi\over 2}, \pi\right)$ on the second $\rp^1$. The element $T\in\pgl(2; \real)$ acts on the second $\rp^1$ of the boundary stratum $(X_1, x, X_5) \cup (x, X_2, X_3, X_4) $ to give $(X_1, x, X_5) \cup (X_4, X_3, X_2, x)$ which is mapped to $(X_1, x, X_5) \cup (x, X_4, X_3, X_2)$ under a rotation in $\pgl(2; \real)$ of the angle between $x$ and $X_4$. The face $(X_1, x, X_5) \cup (x, X_4, X_3, X_2)$ is in the closure of the cell $(X_1, X_4, X_3, X_2, X_5)$ as the part where $X_5$ collides with $X_1$.  In the cell $(X_1, X_4, X_3, X_2, X_5)$, the points $X_1, X_4, X_3$ are fixed at $0,  {\pi \over 2}, \pi$, and the points $(X_2, X_5)$ move in $(\pi, 2\pi) $ with $t_2 \in (\pi, t_5)$ and $t_5 \in (t_2, 2\pi)$. Consequently, $X_5$ can collide with $X_1$ or $X_2$ to form the boundary stratum $(X_1, x, X_5) \cup (x, X_4, X_3, X_2)$ or $(X_2, x, X_5) \cup (x, X_1, X_4, X_3)$, respectively. However, $X_5$ cannot collide with $X_4$ or $X_3$, because $X_4$ and $X_3$ are fixed at ${\pi\over 2}$ and $\pi$, and $X_5$ moves in $(t_2, 2\pi)$ with $t_2\in(\pi,t_5)$. Therefore, the face $(X_2, x, X_5) \cup (x, X_1, X_4, X_3)$ and $(X_1, x, X_5) \cup (x, X_4, X_3, X_2)$ meet at the cell  $(X_1, X_4, X_3, X_2, X_5)$. On the other hand, the face $(X_2, x, X_5) \cup (x, X_1, X_4, X_3)$ is where cells $(X_1, X_4, X_3, X_2, X_5)$ and $(X_1, X_4, X_3, X_5, X_2)$ meet, and the face $(X_1, x, X_5) \cup (x, X_4, X_3, X_2)$ is where cells $(X_1, X_4, X_3, X_2, X_5)$ and $(X_1, X_2, X_3, X_4, X_5)$ meet.
 
\subsection{Categorical Equivalence of Phylogenetic Trees and Stable Curves in $\overline{\mathfrak{M}}_{0,n}$}\label{sec:tree-curve}

The collection of trees $T^\mathcal{X}_x \in\bhv_n$ forms a category. We first present a definition of category attributed to J.~P.~May \cite{may1999concise}, followed by a proof that $T^\mathcal{X}_x \in\bhv_n$ forms a category in Proposition \ref{thm:cat_tree}.

\begin{definition}\rm
A \it category \rm $\mathcal{C}$ consists of a collection of objects, a set $\mathcal{C}(A,B)$ of morphisms (also called maps) between any two objects, an identity morphism $\id_A\in \mathcal{C}(A,A)$ for each object $A$ (usually abbreviated $\id$), and a composition law
$$\circ: \mathcal{C} ( B , C ) \times \mathcal{C} ( A , B )\longrightarrow\mathcal{C} ( A , C )$$
for each triple of objects $A$, $B$, $C$. Composition must be associative, and identity
morphisms must behave as their names dictate: 
$$h\circ (g\circ f)=(h\circ g)\circ f, \ \id \circ f=f, \text{ and } f\circ \id=f$$
whenever the specified composites are defined \cite{may1999concise}.
\end{definition}

The category of trees can be equipped with a few morphisms. The most notable morphism between the objects of trees is the tree isomorphism, which is directly induced from the standard graph isomorphism when the tree is viewed as a graph. 
Two graphs $G_1$ and $G_2$ are isomorphic if there exists a match between their vertices so that two vertices are connected by an edge in $G_1$ if and only if corresponding vertices are connected by an edge in $G_2$. We do not adopt this set of morphisms for our category of trees because each tree in $\bhv_n$ represents the equivalence class of isomorphic trees.

Moerdijk and Weiss have described a category $\Omega$ whose objects are trees of a different type, which allow some edges to have a vertex only on one side \cite{moerdijk2007dendroidal}. The morphisms in their category are defined based on the notion of operad introduced by J.~P.~May \cite{may1997definitions}. The notion of operad encoded how to graft, in which vertices and edges of two trees are placed to obtain a grafted tree combining both trees. The authors viewed any tree $T$ as generating an operad $\Omega(T)$, and the morphisms are operad maps $\Omega(T) \to \Omega(T^\prime)$. We do not adopt this set of morphisms because their trees permit edges to have a vertex only on one side. We define a morphism between two trees if the set of bipartitions of the first tree contains the second. Now we present the notion of homotopy and a special case of homotopy, called a \it deformation retract. \rm
\begin{definition}\rm
A \it homotopy \rm $h : p \simeq q$ between maps $p,q : X \longrightarrow Y$ is a continuous map $h:X\times I \longrightarrow Y$ such that $h(x,0)=p(x)$ and $h(x,1)=q(x)$, where $I$ is the unit interval [0, 1] \cite{may1999concise}. 
\end{definition}

\begin{definition}\rm
A subspace $A$ of a space $X$ is a \it deformation retract \rm if there is a homotopy $h: X \times I \longrightarrow X$ such that $h ( x , 0 ) = x$, $h ( a , t ) = a$, and $h ( x , 1 ) \in A$ for all $ x \in X$, $a \in A$, and $t \in I$. Such a homotopy is called a deformation of $X$ onto $A$ \cite{may1999concise}.
\end{definition}

\begin{proposition}\rm\label{prop:retract}
Denote by $\mathfrak{W}^x_n$ the moduli space of the set of internal edge length of a phylogenetic tree $T^\mathcal{X}_x$, and
$\tilde{\mathfrak{W}}^x_n$ a subspace of $\mathfrak{W}^x_n$ with the $i$-th coordinate being zero. The edge-shrinking map $F_i$ on the split
$S_n^{A_i|B_i}$ 
is defined as: 
\begin{align*}
F_i: \mathfrak{W}^x_n \times [0,1] & \longrightarrow \tilde{\mathfrak{W}}^x_n\\
\left( \left\{w^{A_1\mid B_1}_{1}, \ldots, w^{A_i\mid B_i}_{i}, \ldots, w^{A_k\mid B_k}_{k}\right\}, t\right) & \mapsto \left\{w^{A_1\mid B_1}_{1}, \ldots, (1- t) w^{A_i\mid B_i}_{i}, \ldots, w^{A_k\mid B_k}_{k}\right\}.
\end{align*}
 Then $F_i$ is a homotopy, and $\mathfrak{W}^x_n$ is a deformation retract onto $ \tilde{\mathfrak{W}}^x_n$.
 \end{proposition}
 \begin{proof}
We define the identity map 
\begin{align*}
p_i: \mathfrak{W}^x_n & \longrightarrow \tilde{\mathfrak{W}}^x_n\\
\left\{w^{A_1\mid B_1}_{1}, \ldots, w^{A_i\mid B_i}_{i}, \ldots, w^{A_k\mid B_k}_{k}\right\} & \mapsto \left\{w^{A_1\mid B_1}_{1}, \ldots, w^{A_i\mid B_i}_{i}, \ldots, w^{A_k\mid B_k}_{k}\right\}
 \end{align*}and the map
\begin{align*}
q_i: \mathfrak{W}^x_n & \longrightarrow \tilde{\mathfrak{W}}^x_n\\
\left\{w^{A_1\mid B_1}_{1}, \ldots, w^{A_i\mid B_i}_{i}, \ldots, w^{A_k\mid B_k}_{k}\right\} & \mapsto \left\{w^{A_1\mid B_1}_{1}, \ldots, w^{A_{i-1}\mid B_{i-1}}_{i}, 0, w^{A_{i+1}\mid B_{i+1}}_{i}, \ldots, w^{A_k\mid B_k}_{k}\right\}.
 \end{align*}
Recall that the set of edge lengths of a phylogenetic tree $T^\mathcal{X}_x$ is denoted by $$\mathcal{W}^\mathcal{X}_x= \left\{w^{A_1\mid B_1}_{1}, \ldots, w^{A_i\mid B_i}_{i}, \ldots, w^{A_k\mid B_k}_{k}\right\}.$$
We denote the set of edge lengths of $T^\mathcal{X}_x$ with the $i$-th coordinate being zero
 by $$\tilde{\mathcal{W}}^\mathcal{X}_x= \left\{w^{A_1\mid B_1}_{1}, \ldots, 0, \ldots, w^{A_k\mid B_k}_{k}\right\}.$$Then $F_i$ is a homotopy $p_i \simeq q_i$ because $F_i$ is continuous, $F_i(\mathcal{W}^\mathcal{X}_x,0)=p_i(\mathcal{W}^\mathcal{X}_x)$ and $F_i(\mathcal{W}^\mathcal{X}_x,1)=q_i(\mathcal{W}^\mathcal{X}_x)$. For any set of edge lengths
 $\tilde{\mathcal{W}}^\mathcal{X}_x \in \tilde{\mathfrak{W}}^x_n$, we have 
$F_i (\tilde{\mathcal{W}}^x_n , t ) = \tilde{\mathcal{W}}^x_n$ for all $t \in I$. Thus, $\mathfrak{W}^x_n$ is a deformation retract onto $ \tilde{\mathfrak{W}}^x_n$.
 \end{proof}
 
 The composition of $F_i$s is still a deformation retract, stated as the following corollary:
\begin{corollary}\rm\label{cor:retract}
Consider two phylogenetic trees $T^\mathcal{X}_x$ and $T^\mathcal{X}_y$ with sets of splits
$\mathcal{S}^\mathcal{X}_x$ and $\mathcal{S}^\mathcal{X}_y$ such that $\mathcal{S}^\mathcal{X}_y\subset \mathcal{S}^\mathcal{X}_x$. We let the index set $\mathcal{I} = \{i_1, \ldots, i_r\}$ denote the indices of the set $\mathcal{S}^\mathcal{X}_x \setminus \mathcal{S}^\mathcal{X}_y = \left\{S_n^{A_{i_1}|B_{i_1}}, \ldots, S_n^{A_{i_r}|B_{i_r}}\right\}$, and we define $F_{\mathcal{I}} = F_{i_1}\circ \cdots \circ F_{i_r}$ as the following:
\begin{align*}
F_{\mathcal{I}}: \mathfrak{W}^x_n &\longrightarrow \tilde{\mathfrak{W}}^x_n\\
T^\mathcal{X}_x & \mapsto F_{i_1} \circ \cdots F_{i_r}(T^\mathcal{X}_x),
\end{align*}
where $\mathfrak{W}^x_n$ is the moduli space of the set of internal edge lengths of the tree topology $T^\mathcal{X}_x$, and
$\tilde{\mathfrak{W}}^x_n$ is a subspace of $\mathfrak{W}^x_n$ with the $i_1, \ldots, i_r$-th coordinates being zero.
Then $F_\mathcal{I}$ is a homotopy, and $\mathfrak{W}^x_n$ is a deformation retract onto $ \tilde{\mathfrak{W}}^x_n$.\end{corollary}

We define the category of trees $\mathfrak{T}_n$ with objects $\mathcal{O}(\mathfrak{T}_n)$ being the collection of trees $T^\mathcal{X}_x(\mathcal{S}^\mathcal{X}_x, \mathcal{W}^\mathcal{X}_x)$ in $\bhv_n$. For two objects $T^\mathcal{X}_x(\mathcal{S}^\mathcal{X}_x, \mathcal{W}^\mathcal{X}_x), T^\mathcal{X}_y(\mathcal{S}^\mathcal{X}_y, \mathcal{W}^\mathcal{X}_y) \in \mathcal{O}(\mathfrak{T}_n)$, we define a morphism $\mor(T^\mathcal{X}_x, T^\mathcal{X}_y)$ between $T^\mathcal{X}_x$ and $T^\mathcal{X}_y$ if $\mathcal{S}^\mathcal{X}_y \subset \mathcal{S}^\mathcal{X}_x$.
If $\mathcal{S}^\mathcal{X}_x = \mathcal{S}^\mathcal{X}_y$, we define the morphism $\mor(T^\mathcal{X}_x, T^\mathcal{X}_y)$ being the identity.

\begin{proposition}\rm\label{thm:cat_tree}
$\mathfrak{T}_n$ is a category.
\end{proposition}
\begin{proof}
By definition, $\mor(T^\mathcal{X}_x, T^\mathcal{X}_x)$ is the identity. Therefore, we have an identity morphism $\id_{T^\mathcal{X}_x}\in \mathcal{M}(\ft_n)$ for each object $T^\mathcal{X}_x$. Then we need to show that we have the composition law:
$$\circ: \mor( T^\mathcal{X}_x, T^\mathcal{X}_y ) \times \mor( T^\mathcal{X}_y, T^\mathcal{X}_z )\longrightarrow\mor( T^\mathcal{X}_x , T^\mathcal{X}_z )$$
for each triple of objects $T^\mathcal{X}_x$, $T^\mathcal{X}_y$, $T^\mathcal{X}_z$. By definition, there is a morphism between $T^\mathcal{X}_x$ and $T^\mathcal{X}_y$ if $\mathcal{S}^\mathcal{X}_y\subset \mathcal{S}^\mathcal{X}_x$ and a morphism between $T^\mathcal{X}_y$ and $T^\mathcal{X}_z$ if $\mathcal{S}^\mathcal{X}_z\subset \mathcal{S}^\mathcal{X}_y$.
Thus, we have $\mathcal{S}^\mathcal{X}_z\subset \mathcal{S}^\mathcal{X}_x$ if both are true, and the composition law holds. The composition is associative since
\begin{equation*}
(\mor ( T^\mathcal{X}_x, T^\mathcal{X}_y ) \circ \mor ( T^\mathcal{X}_y, T^\mathcal{X}_z ))  \circ \mor ( T^\mathcal{X}_z, T^\mathcal{X}_w ) \longrightarrow \mor ( T^\mathcal{X}_x, T^\mathcal{X}_z ) \circ \mor ( T^\mathcal{X}_z, T^\mathcal{X}_w ) \longrightarrow
\mor ( T^\mathcal{X}_x, T^\mathcal{X}_w ).
\end{equation*}
We also have
\begin{equation*}
\mor ( T^\mathcal{X}_x, T^\mathcal{X}_y ) \circ (\mor ( T^\mathcal{X}_y, T^\mathcal{X}_z ) \circ \mor ( T^\mathcal{X}_z, T^\mathcal{X}_w ) )\longrightarrow\mor ( T^\mathcal{X}_x, T^\mathcal{X}_y ) \circ \mor ( T^\mathcal{X}_y, T^\mathcal{X}_w )\longrightarrow
\mor ( T^\mathcal{X}_x, T^\mathcal{X}_w ).
\end{equation*}
Thus, we conclude that the composition is associative.
The identity
$$\mor ( T^\mathcal{X}_x, T^\mathcal{X}_x ) \circ \mor ( T^\mathcal{X}_x, T^\mathcal{X}_y ) =\mor ( T^\mathcal{X}_x, T^\mathcal{X}_y )$$
and
$$\mor ( T^\mathcal{X}_x, T^\mathcal{X}_y ) \circ \mor ( T^\mathcal{X}_x, T^\mathcal{X}_x ) =\mor ( T^\mathcal{X}_x, T^\mathcal{X}_y )$$
also follows.
\end{proof}

We transfer our attention to stable curves in $\partial\overline{\mathfrak{M}}_{0,n}(\mathbb{C})$ for now and progress to the connection between the category of trees and the category of stable curves. 
According to Definition \ref{def:nodal}, a nodal curve is a complete algebraic curve such that every one of its points is either smooth or is locally complex-analytically isomorphic to a neighborhood of the origin in the locus with equation $xy = 0$ in $\mathbb{C}^2$, and a stable curve is nodal by Definition \ref{def:stable}.

We denote the category of stable curves with $n$ marked points by $\fc_n$, and the set of objects $\mathcal{O}(\fc_n)$ being the set of stable curves with $n$ marked points defined as the following:

\begin{definition} \normalfont \label{defn:curve}
We define $\mathcal{O}(\fc_n)$ the set of stable curves with $n$ marked points, with elements in $\partial\overline{\mathfrak{M}}_{0,n}(\mathbb{C})=\overline{\mathfrak{M}}_{0,n}(\mathbb{C}) - \mathfrak{M}_{0,n}(\mathbb{C})$. A curve $\mathcal{C}_x^X(\mathcal{S}^X_x)$ in $\mathcal{O}(\fc_n)$ with $k$ components has $k-1$ intersections because stable curves are nodal. Each intersection separates the set of marked points $X$ into two sets, $A_i$ and $B_i$, such that $A_i \cup B_i = X$ and $A_i \cap B_i = \emptyset$. Also, we have $|A_i| \geq 2$ and $|B_i|\geq 2$  under the stable condition. We denote the set of bipartitions induced by intersections of components by
$$\mathcal{S}^X_x = \left\{S^{A_1\mid B_1}_{n}, \ldots, S^{A_k\mid B_k}_{n}\right\},$$
 where $k$ is the number of intersections of the stable curve $\mathcal{C}_x^X(\mathcal{S}^X_x)$, and $x$ labels different stable curves. 
\end{definition}

\begin{figure}
 \begin{center}
  \includegraphics[width=.9\textwidth]{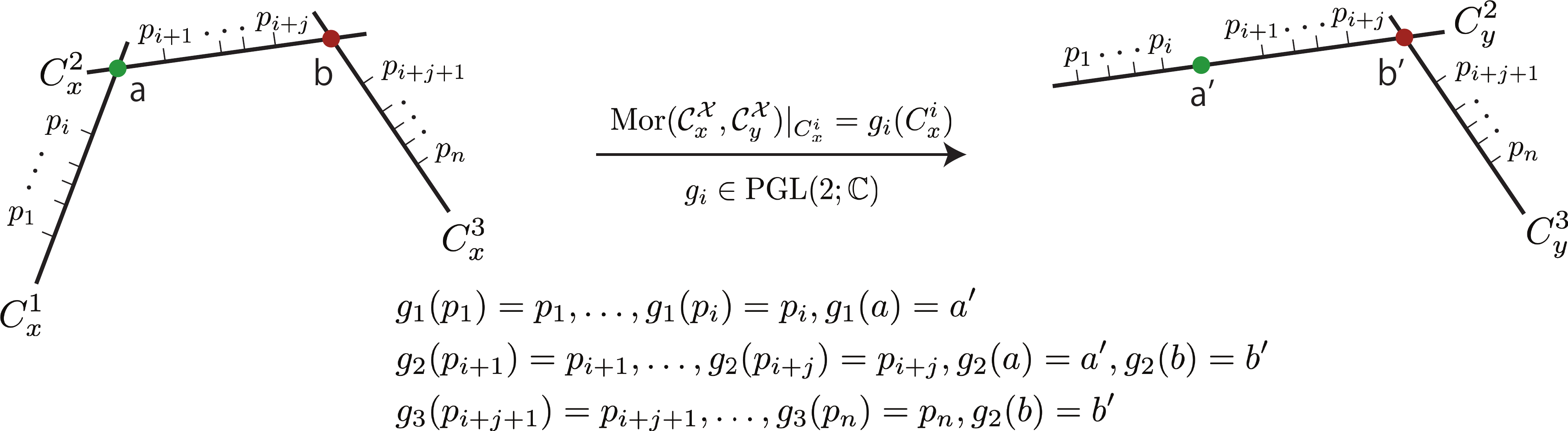}
 \end{center}
 \caption{Morphisms in the category of genus 0 stable curves with $n$ marked points. \it Left: \rm a stable curve $\mathcal{C}_x^X$ with three $\mathbb{P}^1$ components; \it right: \rm a stable curve $\mathcal{C}_y^X$ with two $\mathbb{P}^1$ components.\rm }\label{fig:mor}
\end{figure}

We define a morphism between two objects $\mathcal{C}_x^X, \mathcal{C}_y^X \in \mathcal{O}(\fc_n)$ if for each component $C^i_x$ in $\mathcal{C}_x^X$, there exists $g_i \in \pgl(2; \complex)$ that maps all the marked points on this component to one single component of $\mathcal{C}_y^X$. Also, each intersection of two components $C^i_x$ and $C^j_x$ of $\mathcal{C}_x^X$ is mapped to the same point in $\mathcal{C}_y^X$ under $g_i$ and $g_j$, respectively. The morphism $\mor(\mathcal{C}_x^X, \mathcal{C}_y^X)$ between two stable curves $\mathcal{C}_x^X$ and $\mathcal{C}_y^X$ are these $g_i$s glued along the intersection points:
\begin{align*}
\mor(\mathcal{C}_x^X, \mathcal{C}_y^X): \mathcal{C}_x^X &\longrightarrow \mathcal{C}_y^X\\
C^i_x & \mapsto g_i(C^i_x)
\end{align*}
with
\begin{align*}
g_i : C^i_x \subset \mathcal{C}_x^X&\longrightarrow C^i_y\subset \mathcal{C}_y^X\\
g_i(p_j) & \mapsto p_j\\
g_i(a) &\mapsto a^\prime
\end{align*}
for any marked point $p_j$, the intersection point $a$ of two components $C_x^i, C_x^j\in\mathcal{C}_x^X$, and the intersection point $a^\prime$ of two components $C_y^i, C_y^j\in\mathcal{C}_y^X$ such that $g_i(a) = a^\prime, g_j(a) = a^\prime$. An example is shown as Figure \ref{fig:mor}.

\begin{proposition}\rm\label{thm:cat_curve}
$\fc_n$ is a category.
\end{proposition}
\begin{proof}
First we show that we have an identity morphism $\id_{\mathcal{C}_x^X}\in \mathcal{M}(\fc_n)$ for each object $\mathcal{C}_x^X\in\mathcal{O}(\fc_n)$. This is true because we can set each $g_i$ as the identity in $\pgl(2;\complex)$, so
$$\mor ( \mathcal{C}_x^X, \mathcal{C}_x^X ) \circ \mor ( \mathcal{C}_x^X, \mathcal{C}_y^X ) =\mor ( \mathcal{C}_x^X, \mathcal{C}_y^X )$$
and $$\mor ( \mathcal{C}_x^X, \mathcal{C}_y^X ) \circ \mor ( \mathcal{C}_x^X, \mathcal{C}_x^X ) =\mor ( \mathcal{C}_x^X, \mathcal{C}_y^X ).$$
We need to show that we have the composition law:
$$\circ: \mor( \mathcal{C}_x^X, \mathcal{C}_y^X ) \times \mor ( \mathcal{C}_y^X, \mathcal{C}_z^X )\longrightarrow\mor (\mathcal{C}_x^X , \mathcal{C}_z^X )$$
for each triple of objects $\mathcal{C}_x^X$, $\mathcal{C}_y^X$, $\mathcal{C}_z^X$. We consider the morphism $\mor(\mathcal{C}_x^X, \mathcal{C}_y^X)$ between two stable curves $\mathcal{C}_x^X$ and $\mathcal{C}_y^X$:\begin{align*}
\mor(\mathcal{C}_x^X, \mathcal{C}_y^X): \mathcal{C}_x^X &\longrightarrow \mathcal{C}_y^X\\
C_x^i & \mapsto g_i(C_x^i)
\end{align*}
with
\begin{align*}
g_i : C_x^i \subset \mathcal{C}_x^X&\longrightarrow C_y^i\subset \mathcal{C}_y^X\\
g_i(p_j) & \mapsto p_j\\
g_i(a) &\mapsto a^\prime
\end{align*}
for $g_i\in \pgl(2; \complex)$, any marked point $p_j$, the intersection point $a$ of two components $C_x^i, C_x^j\in\mathcal{C}_x^X$, and the intersection point $a^\prime$ of two components $C_y^i, C_y^j\in\mathcal{C}_y^X$ such that $g_i(a) = a^\prime, g_j(a) = a^\prime$.

The morphism $\mor(\mathcal{C}_y^X, \mathcal{C}_z^X)$ between two stable curves $\mathcal{C}_y^X$ and $\mathcal{C}_z^X$ is:\begin{align*}
\mor(\mathcal{C}_y^X, \mathcal{C}_z^X): \mathcal{C}_y^X &\longrightarrow \mathcal{C}_z^X\\
C_y^i & \mapsto f_i(C_y^i)
\end{align*}
with
\begin{align*}
f_i : C_y^i\subset \mathcal{C}_y^X &\longrightarrow C_z^i\subset \mathcal{C}_z^X\\
f_i(p_j) & \mapsto p_j\\
f_i(b) &\mapsto b^\prime
\end{align*}
for $f_i\in \pgl(2; \complex)$, any marked point $p_j$, the intersection point $b$ of two components $C_y^i, C_y^j\in\mathcal{C}_y^X$, and the intersection point $b^\prime$ of two components $C_z^i, C_z^j\in\mathcal{C}_z^X$ such that $f_i(b) = b^\prime, f_j(b) = b^\prime$.

We need to check that we have the morphism $\mor(\mathcal{C}_x^X, \mathcal{C}_z^X)$ such that
\begin{align*}
\mor(\mathcal{C}_x^X, \mathcal{C}_z^X): \mathcal{C}_x^X &\longrightarrow \mathcal{C}_z^X\\
C_x^i & \mapsto q_i(C_x^i)
\end{align*}
for some $q_i\in \pgl(2; \complex)$ with the desired properties. We set $q_i = f_i \circ g_i$, which satisfies:
\begin{align*}
q_i : C_x^i\subset \mathcal{C}_x^X &\longrightarrow C_z^i\subset \mathcal{C}_z^X\\
q_i(p_j) & \mapsto p_j\\
q_i(a) &\mapsto b^\prime
\end{align*}due to the composition. Meanwhile, $q_i = f_i\circ g_i\in \pgl(2; \complex)$ because $f_i$ and $g_i$ are elements of the group $ \pgl(2; \complex)$. The intersection of component $C_i\in \mathcal{C}_x^X$ is mapped to a single point $b^\prime$ in $\mathcal{C}_z^X$ by all $q_i$s as well.

The associativity of the composition follows from the associativity of the composition in $\pgl(2; \complex)$. Given $\mor(\mathcal{C}_z^X, \mathcal{C}_w^X)$ with
\begin{align*}
\mor(\mathcal{C}_z^X, \mathcal{C}_w^X): \mathcal{C}_z^X &\longrightarrow \mathcal{C}_w^X\\
C_z^i & \mapsto h_i(C_w^i)
\end{align*}
with
\begin{align*}
h_i : C_z^i \subset \mathcal{C}_z^X&\longrightarrow C_w^i\subset \mathcal{C}_w^X\\
h_i(p_j) & \mapsto p_j\\
h_i(c) &\mapsto c^\prime
\end{align*}
for some $h_i\in \pgl(2; \complex)$, then we have
$$h_i \circ (f_i \circ g_i) = (h_i \circ f_i) \circ g_i.$$

Because $g_i$s, $f_i$s, and $h_i$s all send special points to special points, so does their composition. On the left-hand side, we have $f_i \circ g_i$ mapping marked points in $\mathcal{C}_x^X$ to the same marked points in $\mathcal{C}_z^X$, and mapping intersection points in $\mathcal{C}_x^X$ to intersection points or interior points in $\mathcal{C}_z^X$ that are consistent with $f_j \circ g_j$ of the other component $C^j_x$ of $\mathcal{C}_x^X$ that intersects with $C^i_x$. It is  similar for $h_i$ mapping from $\mathcal{C}_z^X$ to $\mathcal{C}_w^X$.
On the right-hand side, we have $g_i$ mapping marked points in $\mathcal{C}_x^X$ to the same marked points in $\mathcal{C}_y^X$, and mapping intersection points in $\mathcal{C}_x^X$ to intersection points or interior points in $\mathcal{C}_y^X$ consistently, and similarly for $h_i \circ f_i$ mapping from $\mathcal{C}_y^X$ to $\mathcal{C}_w^X$. Therefore, we have
$$
\mor ( \mathcal{C}_x^X, \mathcal{C}_y^X ) \circ (\mor ( \mathcal{C}_y^X, \mathcal{C}_z^X ) \circ \mor ( \mathcal{C}_z^X, \mathcal{C}_w^X ) )
=
(\mor ( \mathcal{C}_x^X, \mathcal{C}_y^X ) \circ \mor ( \mathcal{C}_y^X, \mathcal{C}_z^X ))  \circ \mor ( \mathcal{C}_z^X, \mathcal{C}_w^X ).$$
\end{proof}

The category of trees and the category of stable curves are connected through the space of tree topology and the partition space of stable curves that encode the respective combinatorial properties. We first define these spaces, and we then introduce their functorial relationships. We define the category of tree topologies $\overline{\mathfrak{T}}_n$ with objects $\mathcal{O}(\overline{\mathfrak{T}}_n)$ as the collection of equivalence classes of trees $\overline{T}^\mathcal{X}_x(\mathcal{S}^\mathcal{X}_x)$ in $\bhv_n$ that share the same set of splits $\mathcal{S}^\mathcal{X}_x$ but with different weights. For two objects $\overline{T}^\mathcal{X}_x(\mathcal{S}^\mathcal{X}_x), \overline{T}^\mathcal{X}_y(\mathcal{S}^\mathcal{X}_y) \in \mathcal{O}(\overline{\mathfrak{T}}_n)$, we define a morphism $\mor\left(\overline{T}^\mathcal{X}_x, \overline{T}^\mathcal{X}_y\right)$ between $\overline{T}^\mathcal{X}_x$ and $\overline{T}^\mathcal{X}_y$ if $\mathcal{S}^\mathcal{X}_y\subset \mathcal{S}^\mathcal{X}_x$ as the following: Let $\mathcal{I} = \{i_1, \ldots, i_r\}$ denote the indices of the set $\mathcal{S}^\mathcal{X}_x \setminus \mathcal{S}^\mathcal{X}_y = \left\{S_n^{A_{i_1}|B_{i_1}}, \ldots, S_n^{A_{i_r}|B_{i_r}}\right\}$; we define a projection map $\pi_i$ that forgets the $i$-th component, and we denote $\pi_{\mathcal{I}}$ as the composition of projection maps $\pi_{i_1} \circ \cdots \circ \pi_{i_r}$, then
\begin{align*}
\mor(\overline{T}^\mathcal{X}_x, \overline{T}^\mathcal{X}_y):\mathcal{O}(\ft_n) &\longrightarrow \mathcal{O}(\ft_n)\\
\overline{T}^\mathcal{X}_x(\mathcal{S}^\mathcal{X}_x) & \mapsto \pi_\mathcal{I} \left(\overline{T}^x\left(\mathcal{S}^\mathcal{X}_x\right)\right) = \overline{T}^\mathcal{X}_y(\mathcal{S}^\mathcal{X}_y).
\end{align*}
If $\mathcal{S}^\mathcal{X}_x = \mathcal{S}^\mathcal{X}_y$, we define the morphism $\mor\left(\overline{T}^\mathcal{X}_x, \overline{T}^\mathcal{X}_y\right)$ as the identity. 
The proof of $\overline{\mathfrak{T}}_n$ being a category follows from the proof of $\mathfrak{T}_n$ being a category.

\begin{definition}\rm
A \it functor \rm $F$ from a category $\mathcal{C}$ to a category $\mathcal{D}$ assigns to each object $X$ in $\mathcal{C}$ an object $F(X)$ in $\mathcal{D}$
 and to each morphism $f \in\mor(X,Y)$ in $\mathcal{C}$ a morphism $F(f)\in \mor(􏰁F(X),F(Y)􏰂)$ in $\mathcal{D}$, such that $$F(\id_X)=\id_{F(X)} \text{ and } F(f\circ g)=F(f)\circ F(g).$$ This defines a covariant functor. A contravariant functor would differ from this by assigning to $f \in \mor(X,Y)$ a `backwards' morphism $F(f) \in \mor(F(Y), F(X))$ with the same criteria \cite{allen2002cambridge}.
\end{definition}

\begin{proposition}\rm\label{prop:t_to_c}
The functor 
\begin{align*}
F : \mathfrak{C}_n & \longrightarrow \overline{\mathfrak{T}}_n \\
\mathcal{C}_x^X(\mathcal{S}^X_x)
&\mapsto\overline{T}^\mathcal{X}_x(\mathcal{S}^\mathcal{X}_x) \\
\mor(\mathcal{C}_x^X,\mathcal{C}_y^X)
&\mapsto \mor\left(\overline{T}^\mathcal{X}_x, \overline{T}^\mathcal{X}_y\right)
\end{align*}
maps from the category $\fc_n$ of genus 0 stable curves with $n$ marked points to the category of tree topologies $\overline{\mathfrak{T}}_n$.
\end{proposition}
\begin{proof}
The functor $F$ assigns to each object $\mathcal{C}_x^X(\mathcal{S}^X_x)$ in $\fc_n$ an object $F(\mathcal{C}_x^X)$ of $ \overline{\mathfrak{T}}_n$, which is $\overline{T}^\mathcal{X}_x(\mathcal{S}^\mathcal{X}_x)$. $F$ assigns to each morphism 
$\mor(\mathcal{C}_x^X, \mathcal{C}_y^X) \in \mathcal{M}(\fc_n)$ a morphism
$F(\mor(\mathcal{C}_x^X, \mathcal{C}_y^X) )= \mor\left(\overline{T}^\mathcal{X}_x, \overline{T}^\mathcal{X}_y\right)\in \mor\left(\overline{\mathfrak{T}}_n\right) $, such that
$$F\left(\id_{\mathcal{C}_x^X}\right) 
=
 F(\mor(\mathcal{C}_x^X, \mathcal{C}_x^X))
 =\mor\left(\overline{T}^\mathcal{X}_x, \overline{T}^\mathcal{X}_x\right)
 = \id_{F(\mathcal{C}_x^X)}.$$
Now we consider that two morphisms
$\mor(\mathcal{C}_x^X, \mathcal{C}_y^X), \mor(\mathcal{C}_y^X, \mathcal{C}_z^X) \in \mathcal{M}(\fc_n)$, and we have
\begin{align*}
F\left(\mor(\mathcal{C}_x^X, \mathcal{C}_y^X)\circ \mor(\mathcal{C}_y^X, \mathcal{C}_z^X)\right)
& =F(\mor(\mathcal{C}_x^X, \mathcal{C}_z^X))\\
&= \mor\left(\overline{T}^\mathcal{X}_x, \overline{T}^\mathcal{X}_z\right),\\
 F(\mor(\mathcal{C}_x^X, \mathcal{C}_y^X)) \circ F(\mor(\mathcal{C}_x^X, \mathcal{C}_y^X)) &= \mor\left(\overline{T}^\mathcal{X}_x, \overline{T}^\mathcal{X}_y\right)\circ \mor\left(\overline{T}^\mathcal{X}_y, \overline{T}^\mathcal{X}_z\right)\\
 &= \mor\left(\overline{T}^\mathcal{X}_x, \overline{T}^\mathcal{X}_z\right).
 \end{align*}
 Thus, we conclude that $F$ is a functor from $\mathfrak{C}_n$
 to $\overline{\mathfrak{T}}_n$.

\end{proof}

We define the category of the partition space of genus 0 stable curves with $n$ marked points $\ofc_n$ with objects $\mathcal{O}(\ofc_n)$ as the collection of equivalence classes of genus 0 stable curves $\overline{\mathcal{C}}^X_x(\mathcal{S}^X_x)$ in $\mathfrak{C}_n$ with the set of intersections $\mathcal{S}^\mathcal{X}_x$. 
We define the fusion map $\rho_i$ whose effect is forgetting the intersection $S^{A_i|B_i}_n$ of two components $C^i_x$ and $C^j_x$ of $\mathcal{C}^x_n$:
\begin{align*}
\rho_i: \mathcal{O}(\ofc_n)&\longrightarrow \mathcal{O}(\ofc_n)\\
\overline{\mathcal{C}}^X_x\left(\left\{S^{A_1|B_1}_n, \ldots, S^{A_i|B_i}_n, \ldots, S^{A_k|B_k}_n\right\}\right) & \mapsto\overline{\mathcal{C}}^X_x\left(\left\{S^{A_1|B_1}_n, \ldots, \hat S^{A_i|B_i}_n, \ldots, S^{A_k|B_k}_n\right\}\right) .
\end{align*}
For two objects $\overline{\mathcal{C}}^X_x(\mathcal{S}^X_x), \overline{\mathcal{C}}^X_y(\mathcal{S}^X_y) \in \mathcal{O}(\overline{\mathfrak{C}}_n)$, if $\mathcal{S}^X_y \subset \mathcal{S}^X_x$, we define a morphism $\mor\left(\overline{\mathcal{C}}^X_x, \overline{\mathcal{C}}^X_y\right)$ between $\overline{\mathcal{C}}^X_x$ and $\overline{\mathcal{C}}^X_y$ by:
\begin{align*}
\mor\left(\overline{\mathcal{C}}^X_x, \overline{\mathcal{C}}^X_y\right):\mathcal{O}(\fc_n) &\longrightarrow \mathcal{O}(\fc_n)\\
\overline{\mathcal{C}}^X_x(\mathcal{S}^X_x) & \mapsto \rho_\mathcal{I} \left(\overline{\mathcal{C}}^x(\mathcal{S}^X_x)\right),
\end{align*}
where $\rho_{\mathcal{I}}$ is the composition of fusion maps $\rho_{i_1} \circ \cdots \circ \rho_{i_r}$.
If $\mathcal{S}^X_x = \mathcal{S}^X_y$, we define the morphism $\mor\left(\overline{\mathcal{C}}^X_x, \overline{\mathcal{C}}^X_y\right)$ to be the identity. 

\begin{proposition}\rm\label{prop:c_to_t}
The functor 
\begin{align*}
G : \mathfrak{T}_n & \longrightarrow \overline{\mathfrak{C}}_n \\
T^\mathcal{X}_x(\mathcal{S}^\mathcal{X}_x, \mathcal{W}^\mathcal{X}_x)
&\mapsto\overline{\mathcal{C}}^X_x(\mathcal{S}^X_x) \\
\mor(T^\mathcal{X}_x,T^\mathcal{X}_y)
&\mapsto \mor\left(\overline{\mathcal{C}}^X_x, \overline{\mathcal{C}}^X_y\right)
\end{align*}
maps from the category of trees $\mathfrak{T}_n$
to the category $\ofc_n$ of the partition space of genus 0 stable curves with $n$ marked points.
\end{proposition}
\begin{proof}
The functor $G$ assigns to each object $T^\mathcal{X}_x(\mathcal{S}^\mathcal{X}_x, \mathcal{W}^\mathcal{X}_x)$ in $\ft_n$ an object $$G(T^\mathcal{X}_x)=\overline{\mathcal{C}}^X_x(\mathcal{S}^X_x)\in\ofc_n.$$
$G$ assigns to each morphism 
$\mor(T^\mathcal{X}_x, T^\mathcal{X}_y) \in \mathcal{M}(\ft_n)$ a morphism
$$G(\mor(T^\mathcal{X}_x, T^\mathcal{X}_y) )= \mor\left(\overline{\mathcal{C}}^X_x,\overline{\mathcal{C}}^X_y\right)\in \mathcal{M}\left(\overline{\mathfrak{C}}_n\right),$$
such that
$$G\left(\id_{T^\mathcal{X}_x}\right) =G\left(\mor(T^\mathcal{X}_x, T^\mathcal{X}_x)\right)
=\mor\left( \overline{\mathcal{C}}^X_x, \overline{\mathcal{C}}^X_x\right)
 = \id_{G(T^\mathcal{X}_x)}.$$
Now we consider that two morphisms
$\mor(T^\mathcal{X}_x,T^\mathcal{X}_y), \mor(T^\mathcal{X}_y, T^\mathcal{X}_z) \in\mathcal{M}(\ft_n)$, and similarly we have
$$G\left(\mor(T^\mathcal{X}_x, T^\mathcal{X}_y)\circ \mor(T^\mathcal{X}_y, T^\mathcal{X}_z)\right)
= G(\mor(T^\mathcal{X}_x, T^\mathcal{X}_y)) \circ G(\mor(T^\mathcal{X}_y, T^\mathcal{X}_z)).$$
Thus, we can conclude that $G$ is a functor from $\mathfrak{T}_n$ to $\ofc_n$.
\end{proof}

Our next goal is to prove that the category of the partition space of genus 0 stable curves with $n$ marked points $\overline{\fc}_n$ is equivalent to the category of tree topologies $\overline{\mathfrak{T}}_n$. We first define two functors between these two categories.
The functor 
\begin{align*}
\overline{F} :\overline{\mathfrak{C}}_n & \longrightarrow \overline{\mathfrak{T}}_n \\
\overline{\mathcal{C}}^X_x(\mathcal{S}^X_x)
&\mapsto\overline{T}^\mathcal{X}_x(\mathcal{S}^\mathcal{X}_x) \\
\mor(\overline{\mathcal{C}}^X_x,\overline{\mathcal{C}}^X_y)
&\mapsto \mor\left(\overline{T}^\mathcal{X}_x, \overline{T}^\mathcal{X}_y\right)
\end{align*}
maps from the category  of the partition space of genus 0 stable curves with $n$ marked points $\ofc_n$ to the category of tree topologies $\overline{\mathfrak{T}}_n$.
The functor 
\begin{align*}
\overline{G} : \overline{\mathfrak{T}}_n & \longrightarrow \overline{\mathfrak{C}}_n \\
\overline{T}^\mathcal{X}_x(\mathcal{S}^\mathcal{X}_x)&\mapsto\overline{\mathcal{C}}^X_x(\mathcal{S}^X_x) \\
\mor\left(\overline{T}^\mathcal{X}_x,\overline{T}^\mathcal{X}_y\right)
&\mapsto \mor\left(\overline{\mathcal{C}}^X_x, \overline{\mathcal{C}}^X_y\right)
\end{align*}
maps from the category of tree topologies $\overline{\mathfrak{T}}_n$ to the category of the partition space of genus 0 stable curves with $n$ marked points $\ofc_n$.

\begin{definition}\rm
A \it natural transformation \rm $\alpha: F \longrightarrow G$ between functors $ \mathcal{C}\longrightarrow \mathcal{D}$ is a map of functors. It consists of a morphism $\alpha_A: F(A) \longrightarrow G(A)$ for each object $A$ of $ \mathcal{C}$ such that the following diagram commutes for each morphism $f: A \longrightarrow B$ of $\mathcal{C}$ \cite{may1999concise}:
\[
\begin{tikzcd}[column sep=huge,row sep=huge]
F(A) \arrow[r,"F(f)"] \arrow[d,swap,"\alpha_A"] & F(B) \arrow[d,shift left=.75ex,"\alpha_B"] \\
G(A) \arrow[r, "G(f)"] & G(B).
\end{tikzcd}
\]
\end{definition}

\begin{definition}\rm
Two categories $\mathcal{C}$ and $\mathcal{D}$ are \it equivalent \rm if there are functors $F : \mathcal{C} \longrightarrow \mathcal{D}$ and
$G :\mathcal{D} \longrightarrow \mathcal{C}$ and natural isomorphisms $FG \longrightarrow \Id$ and $GF \longrightarrow \Id$, where the $\Id$
are the respective identity functors \cite{may1999concise}.
\end{definition}

\begin{theorem}\rm\label{thm:equiv}
The category $\overline{\mathfrak{T}}_n$ and the category $\ofc_n$ are equivalent.
\end{theorem}
\begin{proof}
To show that the category $\overline{\mathfrak{T}}_n$ and the category $\ofc_n$ are 
equivalent, we need to show that there are natural isomorphisms $\overline{F}\overline{G} \longrightarrow\Id_{\overline{\mathfrak{T}}_n}$ and
$\overline{G}\overline{F} \longrightarrow\Id_{\overline{\mathfrak{C}}_n}$, where the $\Id_{\overline{\mathfrak{T}}_n}$ and $\Id_{\overline{\mathfrak{C}}_n}$ are the respective identity functors on $\overline{\mathfrak{T}}_n$ and $\overline{\mathfrak{C}}_n$.
We first show that there is a natural isomorphism $\alpha_x$ between functors $\overline{F}\overline{G}$ and $\Id_{\overline{\mathfrak{T}}_n}$.
Recall that each intersection of a stable curve with $n$ marked points separates the set of marked points $X$ into two sets $A_i$ and $B_i$, such that $A_i \cup B_i = X$ and $A_i \cap B_i = \emptyset$ with $|A_i| \geq 2$ and $|B_i|\geq 2$ under stable condition. 
On the other hand,  the internal edges of a phylogenetic tree have the property that each edge separates the set $\mathcal{X}$ as the union of two disjoint subsets with both subsets having at least two leaves.  This concludes that the set of intersections  $\mathcal{S}^X$ of stable curves are in one-to-one correspondence with the set of splits $\mathcal{S}^\mathcal{X}$ of phylogenetic trees for some fixed number $n$ of marked points or taxa with $n\geq 4$. This implies that both $\overline{F}$ and $\overline{G}$ are isomorphisms; thus, the composition $\overline{F}\overline{G}$ is an isomorphism as well. Therefore, for each object $\overline{T}^\mathcal{X}_x\in\mathcal{O}(\overline{\mathfrak{T}}_n)$, we have a morphism $\alpha_x: \overline{F}\overline{G}\left(\overline{T}^\mathcal{X}_x\right) \longrightarrow \Id_{\overline{\mathfrak{T}}_n}\left(\overline{T}^\mathcal{X}_x\right)$.
For each morphism $\mor(\overline{T}^\mathcal{X}_x,\overline{T}^\mathcal{X}_y)\in \mathcal{M}(\overline{\mathfrak{T}}_n)$, 
the diagram 
\[
\begin{tikzcd}[column sep=huge,row sep=huge]
\overline{F}\overline{G}(\overline{T}^\mathcal{X}_x) \arrow[r, "\phi"] \arrow[d,swap,"\alpha_x"] &\overline{F}\overline{G}(\overline{T}^\mathcal{X}_y) \arrow[d,shift left=.75ex,"\alpha_y"] \\
\Id_{\overline{\mathfrak{T}}_n}(\overline{T}^\mathcal{X}_x) \arrow[r,"\psi"] & \Id_{\overline{\mathfrak{T}}_n}(\overline{T}^\mathcal{X}_y)
\end{tikzcd}
\]
with
$$\phi= \overline{F}\overline{G} \left(\mor\left(\overline{T}^\mathcal{X}_x, \overline{T}^\mathcal{X}_y\right)
\right), \ \psi= \Id_{\overline{\mathfrak{T}}_n}\left(\mor\left(\overline{T}^\mathcal{X}_x, \overline{T}^\mathcal{X}_y\right)
\right)$$
commutes, since\begin{align*}
\alpha_y \circ \phi \left(\overline{F}\overline{G} \left(\overline{T}^\mathcal{X}_x\right)\right)
 & = \alpha_y \circ\overline{F}\overline{G} \left(\overline{T}^\mathcal{X}_y\right)\\
&=\Id_{\overline{\mathfrak{T}}_n}\left(\overline{T}^\mathcal{X}_y\right),\\
\psi\circ \alpha_x \left( \overline{F}\overline{G} \left(\overline{T}^\mathcal{X}_x\right)\right)
&=\psi\circ  \Id_{\overline{\mathfrak{T}}_n}\left(\overline{T}^\mathcal{X}_x\right)\\
&=\Id_{\overline{\mathfrak{T}}_n}\left(\overline{T}^\mathcal{X}_y\right).
\end{align*}

Similarly, there is a natural isomorphism
$\overline{G}\overline{F} \longrightarrow\Id_{\overline{\mathfrak{C}}_n}$.
Therefore, the category $\overline{\mathfrak{T}}_n$ and the category $\ofc_n$ are equivalent.
\end{proof}

\subsection{Isomorphism Between the Dual Intersection Complex and $\pbhv_n$}\label{sec:dual_tree}
 
In this section, we discuss the dual relation between the space of stable curves and the space of phylogenetic trees. We introduce the notion of the dual intersection complex and present its connection with the projective Billera-Holmes-Vogtmann space of phylogenetic trees $\pbhv_n$.

\begin{theorem}\rm\label{thm:curve_tree_dual}
The space of the stable curves and tree space are the dual of each other such that a tree with $k$ internal edges and $n$ nodes corresponds to the moduli space with $n-3-k$ complex parameters in $\overline{\mathfrak{M}}_{0,n}(\mathbb{C})$.
\end{theorem}
\begin{proof}
We proceed with the proof by induction.
For the base case, the interior of $\overline{\mathfrak{M}}_{0,n}(\mathbb{C})$ corresponds to a tree with no internal edges. Since the dimension of trees is determined by the number of its internal edges, having no internal edge implies that the dimension is zero, i.e., only a point.
Therefore, a tree with $0$ internal edges and $n$ nodes corresponds to the moduli space with $n-3$ complex parameters in $\overline{\mathfrak{M}}_{0,n}(\mathbb{C})$.

For the inductive hypothesis, we assume that a tree with $k$ internal edges corresponds to $n-3-k$ complex parameters, and we show this is also true for $k + 1$. 
We extend a collection of taxa $A_i$ from its complement $B_i$ in $\mathcal{X}$ using an internal edge, preserving existing bipartitions. The new bipartition corresponds to a suitable collision in $\overline{\mathfrak{M}}_{0,n}(\mathbb{C})$, which creates an intersection that separates the set of marked points into $A_i$ and $B_i$. This is a codimension 1 subspace related to the original space. Thus, this gives a space with $n-4-k$ complex parameters, as desired.
\end{proof}

In particular, the dual space theorem gives a sufficient and necessary condition for a genus 0 nodal curve being stable. 

\begin{proposition}\rm\label{prop:remain}
A stable curve remains stable after removing a marked point if and only if it corresponds to a degenerated tree.
\end{proposition}
\begin{proof}If a stable curve remains stable after removing a marked point, it implies that at least one component of the stable curve has more than three special points. Thus, that component corresponds to a node in a tree with valence greater than 3; therefore, the tree is degenerated. Conversely, 
if a stable curve becomes unstable after removing a marked point, this implies that each component of the stable curve has precisely three special points. Therefore, the tree corresponding to this stable curve has precisely valence 3 for each internal node; thus, it is a binary tree, i.e., the tree is not degenerated.
\end{proof}
 
We consider the dual intersection complex of boundary divisors, the irreducible components of the boundary of the moduli space of stable $n$-pointed genus 0 curves. We associate each boundary divisor $\mathcal{C}_i^\mathcal{X}$ with a codimension 1 subspace with a point $p_i\in\mathbb{R}^n$, and we define it as a 0-simplex of the intersection complex. We connect $p_i$ and $p_j$ to form a 1-simplex if $\mathcal{C}_i^\mathcal{X}$ intersects $\mathcal{C}_j^\mathcal{X}$. Furthermore, if $\mathcal{C}_i^\mathcal{X}\cap\mathcal{C}_j^\mathcal{X}\cap \mathcal{C}_k^\mathcal{X} \neq \emptyset$, then $p_i, p_j$, and $p_k$ form a 2-simplex, and so on. Thus, we have:

\begin{proposition}\rm\label{prop:dual}
The dual intersection complex of boundary divisors of $\overline{\mathfrak{M}}_{0,n}(\mathbb{C})$ as an abstract simplicial complex is isomorphic to $\pbhv_{n-1}$ as an abstract simplicial complex. Consequently, their respective geometric realizations are homeomorphic.
\end{proposition}
\begin{proof}
We proceed with the proof by induction. We first consider the 0-simplex case between these two spaces as the base case.  By Definition \ref{def:boundary}, the boundary $\Delta = \overline{\mathfrak{M}}_{0,n}(\mathbb{C}) - \mathfrak{M}_{0,n}(\mathbb{C})$ is a divisor consisting of genus 0 stable curves with one node; we consider the boundary divisors as the 0-simplices of the dual intersection complex. Each component of a boundary divisor has at least two marked points in order to satisfy the stable condition. Thus, the set of boundary divisors corresponds to the set of bipartitions of $n$ points, with each subset having at least two elements.

For the $\pbhv_{n-1}$ counterpart, according to Definition \ref{def:PBHV}, $\pbhv_{n-1}$ inherits a simplicial structure from $\bhv_{n-1}$, where the $k$-simplices of $\pbhv_{n-1}$ are points projected from points in $\bhv_{n-1}$ with exactly $k+1$ nonzero edges. Thus, the 0-simplices of $\pbhv_{n-1}$ are phylogenetic trees with precisely one nonzero edge, and each side of the edge has at least two leaves by Definition \ref{def:tree}. Therefore, given a fixed integer $n\geq 4$, the 0-simplices of the dual intersection complex of boundary divisors of $\overline{\mathfrak{M}}_{0,n}(\mathbb{C})$ are isomorphic to the 0-simplices of $\pbhv_{n-1}$.

For the induction hypothesis, we assume that the collection of $(k-1)$-simplices
formed by $k$ boundary divisors is isomorphic to the collection of $(k-1)$-dimensional simplices of $\pbhv_{n-1}$ as an abstract simplicial complex. Now we will show this is true for the $k$-simplices formed by $k+1$ boundary divisors that intersect each other.

By Definition \ref{defn:compact}, $\overline{\mathfrak{M}}_{g,n}(\mathbb{C})$ compactifies $\mathfrak{M}_{g,n}(\mathbb{C})$ without ever allowing the points to come together. Each time the points come together in one place, the curve sprouts off one component that is isomorphic to the projective line. The intersection between the new component and the original one
induces a bipartition corresponding to a boundary divisor that intersects the previous ones. 
Because the new boundary divisor corresponds to sprouting a new projective line, the bipartitions of the existing boundary divisors are still intact. This implies that the bipartition of the new boundary divisor is compatible with existing bipartitions. Therefore, the bipartition corresponding to the new boundary divisor and the existing $k$ splits form a $k$-simplex of $\pbhv_{n-1}$. Similarly, 
 the process of $(k-1)$-simplices of $\pbhv_{n-1}$  forming $k$-simplices by adding a compatible split corresponds to the process of extending $(k-1)$-simplices of the dual intersection complex into $k$-simplices  by adding a new boundary divisor that intersects the existing ones.
This completes the inductive step.
Therefore, the dual intersection complex of boundary divisors 
 is isomorphic to $\pbhv_{n-1}$ as an abstract simplicial complex. Consequently, their respective geometric realizations are homeomorphic.

\end{proof}

\section{Algebraic Morphisms Between Phylogenetic Networks and Stable Curves in $\overline{\mathfrak{M}}_{g,n}(\mathbb{C})$}\label{sec:network}

Phylogenetic networks, broadly defined, refer to any graph used to represent evolutionary relationships among a set of taxa \cite{huson2010phylogenetic}. Usually, phylogenetic networks are defined as unrooted graphs, such that every split is represented by an array of parallel edges with the same length to model evolutionary phenomena \cite{huson2010phylogenetic}. The internal nodes of a phylogenetic network represent ancestral species, and nodes with more than two parents correspond to reticulate events such as hybridization or recombination \cite{huson2006application}.  Classical phylogenetic methods construct phylogenetic trees and networks from dissimilarity maps, i.e., the distance matrices of taxa. The main methods are the neighbor-joining algorithm to construct phylogenetic trees \cite{saitou1987neighbor} and the neighbor-net algorithm to construct split networks \cite{bryant2002neighbornet}. In our context, phylogenetic networks refer to the split networks constructed by the neighbor-net algorithm, which are also called the circular split network (CSN), to reflect the fact that the bipartitions of phylogenetic networks are compatible with certain circular ordering \cite{devadoss2017space} defined as the following:

\begin{definition}\label{def:circular}\rm
 A set of bipartitions $\mathcal{S}^\mathcal{X}$ is called \it compatible \rm with a cyclic ordering of $\mathcal{X}$ if every split is of the form $S_k^{A_k|B_k}$ with $A_k = \{V_i, V_{i+1}, \ldots, V_j\}, B_k = \{V_{j+1}, \ldots, V_{i-1}\}$ such that $V_{n+1} = V_1$.
A set of bipartitions  $\mathcal{S}^\mathcal{X}$  is called \it circular \rm if it is compatible with some cyclic ordering.
\end{definition}

If the set of bipartition $\mathcal{S}^\mathcal{X}$ is pairwise compatible, then the graph realization of $\mathcal{S}^\mathcal{X}$ is a tree; if $\mathcal{S}^\mathcal{X}$  is compatible with a cyclic ordering, then its graph realization is a network. In particular, a set of pairwise compatible splits can always be arranged to be compatible with a cyclic ordering, i.e., the bipartitions of a phylogenetic tree can always be arranged to be compatible with certain circular ordering.

The phylogenetic networks we study can be viewed as graph realizations of a collection of splits, which are bipartitions of the taxa $\mathcal{X} = A \cup B$ such that $A \cap B = \emptyset$, $|A| \geq 2$, and $|B|\geq 2$. Given a collection of splits, the phylogenetic network representing this set of splits can be generated by the Circular Network Algorithm introduced by Huson \cite{huson2010phylogenetic}, as discussed in Section \ref{subsec:alg}. A phylogenetic network contains the same information as a list of splits such that each split $S^{A_i|B_i}_n$ is equipped with weight $w^{A_i|B_i}_i\in \real_{\geq 0}$ \cite{huson2006application}, and the space of phylogenetic networks has interesting, yet natural, topology from gluing together network spaces that share the same splits. The moduli space of phylogenetic networks is a cubical complex \cite{wu2019comparison}; in this section, we show that the collection of phylogenetic networks forms a category, we discuss the categorical relationships between the space of phylogenetic networks and the space of genus 0 real and complex stable curves in $\overline{\mathfrak{M}}_{0,n}(\mathbb{C})$ and  $\overline{\mathfrak{M}}_{0,n}(\mathbb{R})$, and we prove the space of network topologies forms a category that maps injectively into the partition space of complex high genus stable curves in $\overline{\mathfrak{M}}_{g,n}(\mathbb{C})$.

\subsection{The Space of Circular Split Networks}\label{subsec:networks}

In this section, we provide formal definitions of phylogenetic networks, the space of phylogenetic networks $\csn_n$, and the projectivized space $\pcsn_n$. We then introduce the split network constructing algorithm, along with a handful of examples.

\begin{definition}\rm\label{def:network}
A phylogenetic network (denoted by $N^{\mathcal{X}}$) is a planar graph with extra data that is described momentarily.  With regard to the notation, $\mathcal{X}$ here denotes a set of taxa $\mathcal{X}$ (which is a finite set labeled consecutively by the integers starting from 1). 
The underlying graph has some number of internal edges (edges with two vertices) and $n$ external edges (edges with just a single vertex).  The external edges are in 1-1 correspondence with the set $\mathcal{X}$, and they are called the \it leaves of the graph. \rm Meanwhile, each internal edge has an associated \it length, \rm which is a positive real number.  The set of lengths of the internal edges is denoted by $\mathcal{W}^\mathcal{X}$.  The graph with the internal edge weight assignment is constrained by the requirement that its set of internal edges can be decomposed as a union of disjoint subsets of edges (which are called \it splits\rm) with the following five properties: 
\begin{enumerate}
\item	The edges in any given subset have the same weight. 
\item	Each  edge is the only edge between its end vertices. 
\item	No two edges share an end vertex. 
\item	The end vertices of the edges are on the same two disjoint edge paths, which must be edge paths in the graph that  separates the set $\mathcal{X}$ as the union of the same two disjoint subsets with both subsets having at least two leaves. 
\item	Removing these edges disconnects the graph. 
\end{enumerate}
Edges in the same split are called \it parallel \rm edges.  The set of splits is denoted by $S^\mathcal{X}$; it is an ordered set with components $\left\{S_1^{A_1| B_1}, \ldots, S_k^{A_k | B_k}\right\}$ with any given pair $(A_i, B_i)$ being the corresponding disjoint subset decomposition of $\mathcal{X}$.  (As explained below, the split is a priori determined by the decomposition of $\mathcal{X}$.)  The set $\mathcal{W}^\mathcal{X}$ of weights of the internal edges is sometimes written as $\left\{w_1^{A_1| B_1}, \ldots, w_k^{A_k | B_k}\right\}$ with each constituent being the weight assigned to the edges in the like-labeled split. By way of additional notation:   when two distinct networks have the same set $\mathcal{X}$, we denote them as $N^\mathcal{X}_x$ and $N^\mathcal{X}_y$. 
\end{definition}

The \it space of circular split networks, \rm denoted as $\csn_n$, is the space of isometry classes of split networks with $n$-labeled leaves, where the nonzero weights are on the internal branches, and splits share the same circular ordering. The building blocks of $\csn_n$ are the bipartitions $S_n^{A_i|B_i}$, which can be thought of as a \it basis. \rm The number of bipartitions of $\csn_n$ is:
  \begin{equation*}
    \begin{cases}
      \displaystyle\sum_{i=2}^k {n\choose i} & \text{if $n = 2k+1$ for some integer $k\geq 2$},\\
  \displaystyle\sum_{i=2}^{k-1} {n\choose i} + \left.{n\choose k}\right/2& \text{if $n = 2k$ for some integer $k\geq 2$}.\\
    \end{cases}
\end{equation*}
This combinatorial formula can be simplified as $2^{n-1}-n-1$. 

The number of circular ordering for $n$ labels is ${(n-1)!\over 2}$, which can be regarded as the quotient of $n!$ by the order of the dihedral group $\dih_{2n}$. The space $\csn_n$ is constructed by ${(n-1)!\over 2}$ positive orthants glued along the common bipartitions they share. Since fixing the circular ordering is equivalent to labeling an $n$-gon with the set $\{1, \ldots, n\}$, there are ${n(n-3)\over 2}$ bipartitions compatible with each circular order,  which is the number of chords of an $n$-gon. These ${n(n-3)\over 2}$ bipartitions span an $\mathbb{R}^{n(n-3)\over 2}_{\geq 0}$ space, called an orthant. Phylogenetic networks residing in the interior of this orthant share the same set of splits with positive edge lengths, and they are said to have the same \it network topology. \rm The boundary of orthants in the space $\csn_n$ consists of phylogenetic networks with edges collapsed to length zero, so the set of splits of phylogenetic networks on the boundary of an orthant is a subset of the set of splits of phylogenetic networks on the interior of this orthant. Furthermore, the boundary of an orthant can be viewed as the intersection of orthants, and phylogenetic networks on the boundary can be obtained by collapsing edges from either orthant with this boundary. Thus, the $\csn_n$ space is a cubical complex built by tiling each orthant with unit cubes of dimension ${n(n-3)\over 2}$ in this way \cite{wu2019comparison}.
\begin{figure}[!htb]\centering
 \includegraphics[width=.85\textwidth]{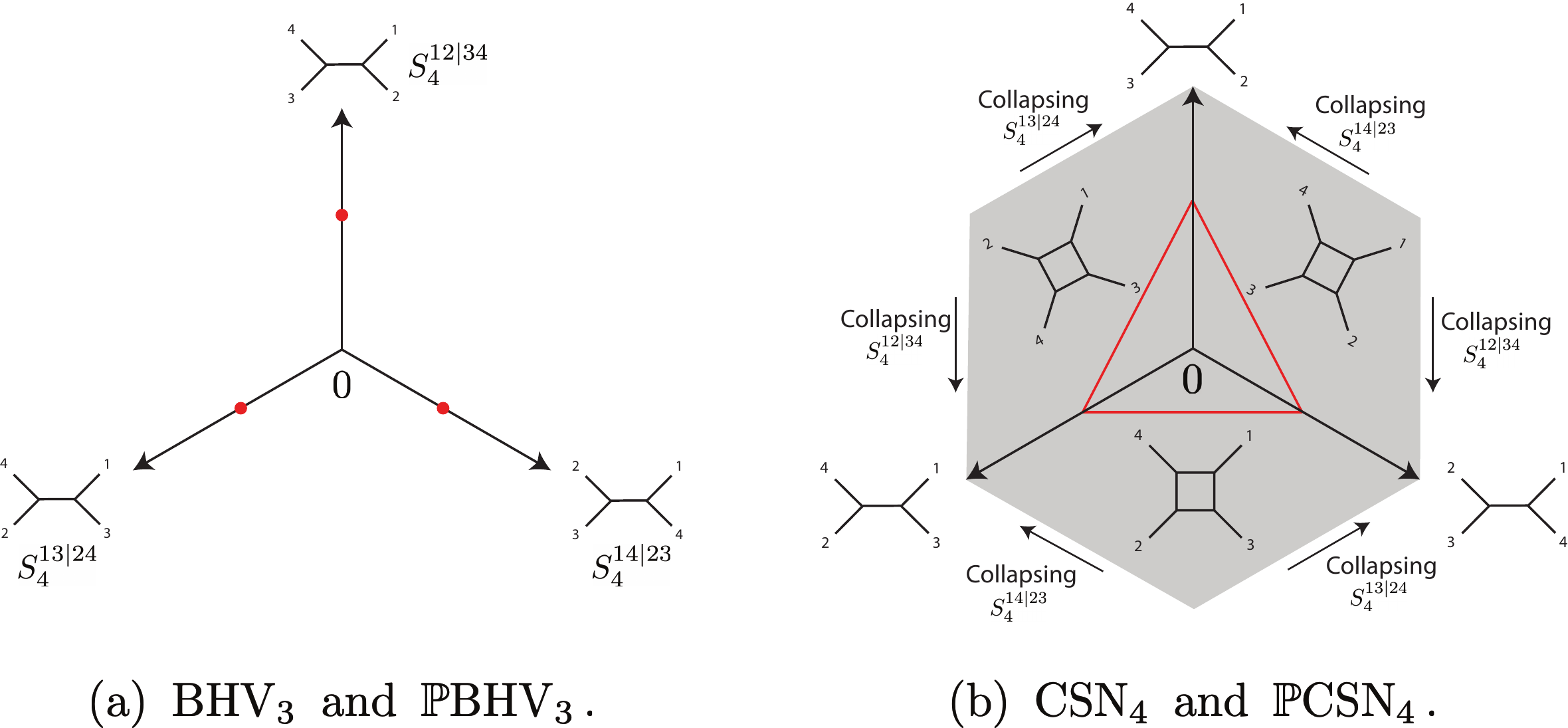}
\caption{Inclusion of $\bhv_3 \subset \csn_4$ and $\pcsn_4 \subset \pbhv_3$. (a) $\bhv_3$: Each ray is a 1-dimensional orthant in $\bhv_3$. The three red dots represent $\pbhv_3$. (b) $\csn_4$: Each ray represents a split; each orthant corresponds to a distinctive network topology. The red triangle represents $\pcsn_4$.}\label{fig:PBHVPCSN}
\end{figure}

The \it projective space of circular split networks $\pcsn_n$ \rm is the co-dimension one subspace of $\csn_n$ such that the sum of the internal edges is 1:
$$\pcsn_n = \left\{
N^\mathcal{X}_x(\mathcal{S}^\mathcal{X}_x,\mathcal{W}^\mathcal{X}_x) \left| \sum_{i=1}^k w_i^{A_i|B_i} = 1
\right\}\right.\subset\csn_n.$$
Each orthant of $\csn_n$ is homeomorphic to $\mathbb{R}^{n(n-3)\over 2}$, and each orthant of $\pcsn_n$ is homeomorphic to $\mathbb{R}^{{n(n-3)\over 2}-1}$. In Figure \ref{fig:PBHVPCSN}a, each ray is a 1-dimensional orthant in $\bhv_3$, and each orthant represents a different tree topology, which corresponds to a distinctive non-trivial split: $S_4^{12|34}$,  $S_4^{13|24}$, and  $S_4^{14|23}$, respectively. $\pbhv_3$ consists of the red dots in Figure \ref{fig:PBHVPCSN}a, with each dot of distance 1 from the origin.
In Figure \ref{fig:PBHVPCSN}b, each ray can be identified with a ray in $\bhv_3$. In general, $\bhv_n$ is a subspace of $\csn_n$, and $\pbhv_n$ is a subspace of $\pcsn_n$ \cite{wu2019comparison}. Each orthant in $\csn_4$ is 2-dimensional, spanned by two rays that share the same circular ordering. $\pcsn_4$ is the subspace of $\csn_4$ marked by the red triangle in the center of Figure \ref{fig:PBHVPCSN}b; for each point on $\pcsn_4$, their coordinates sum up to 1.

$\pcsn_n$ inherits a simplicial structure from $\csn_n$, where the geometric realization of $k$-simplices of $\pcsn_n$ are points projected from points in $\csn_n$ with exactly $k+1$ nonzero edges under the constraint that these edges sum up to 1. Each codimension $m$ face of a $k$-simplex of $\pcsn_n$ is a $(k-m)$-simplex projected from a $(k-m+1)$-simplex of $\csn_n$ spanned by the split set $\mathcal{S}^\mathcal{X}_x \setminus \left\{S_1^{A_1|B_1}, \ldots, S_m^{A_m|B_m}\right\}$ for some set of $m$ splits in $\mathcal{S}^\mathcal{X}_x$.

Furthermore, consider two simplices of $\pcsn_n$ spanned by the set of splits $\mathcal{S}^\mathcal{X}_x$ and $\mathcal{S}^\mathcal{X}_y$; the intersection of these two simplices is spanned by $\mathcal{S}^\mathcal{X}_x \cap \mathcal{S}^\mathcal{X}_y$, which is also a face of each simplex. Each edge of the red triangle in Figure \ref{fig:PBHVPCSN}b is a $1$-simplex of $\pcsn_4$, and the vertices of the red triangle are the 0-simplices. Taking the left orthant as an example, the edge is a subset of $\csn_4$ that contains all networks satisfying $w_4^{12|34} + w_4^{13|24} = 1$. The left-end vertex is the phylogenetic network with one split $S_4^{13|24}$ satisfying $ w_4^{13|24} = 1$, and the right-end vertex is the phylogenetic network containing one split $S_4^{14|23} = 1$ with weight $ w_4^{14|23} = 1$. This shows that the codimension one face of a 1-simplex of $\pcsn_4$ is also a simplex, which is 0-dimensional. Also, the 0-simplices are intersections of the 1-simplices, and each 0-simplex is a face of the 1-simplices connected to it.

\subsection{Split Network Constructing Algorithm}\label{subsec:alg}
\begin{figure}
\includegraphics[width=\textwidth]{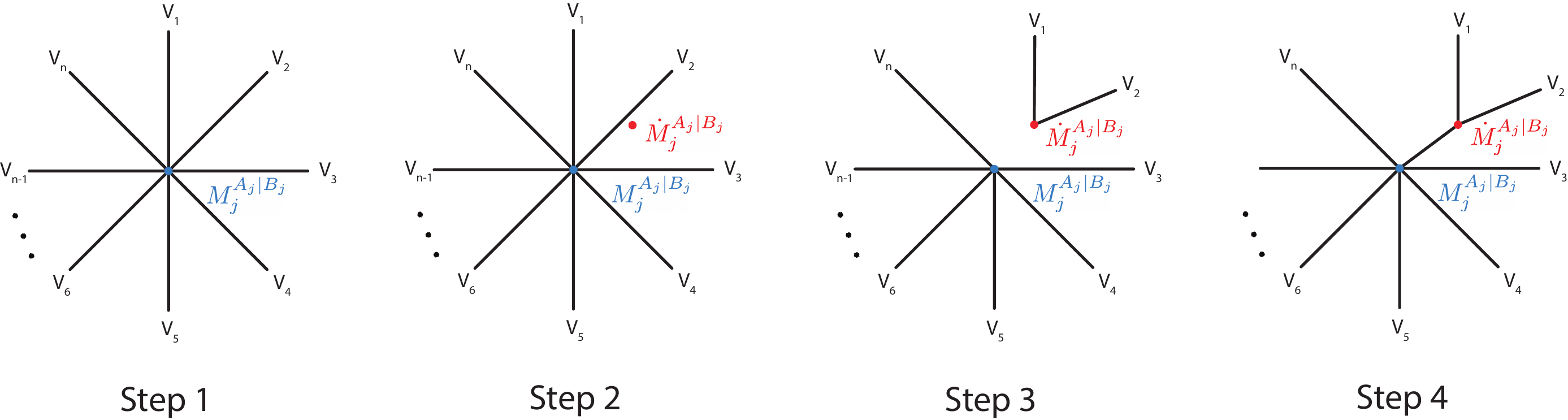}
\caption{Adding the first split for the set of taxa $\mathcal{X} = \{V_1, \ldots, V_n\}$. 
\it Step 1.\rm Determine the split path $M_j^{A_j|B_j}$.
\it Step 2. \rm Create a copy $\dot M_j^{A_j|B_j}$ of $M_j^{A_j|B_j}$.
\it Step 3. \rm Redirect nodes.
\it Step 4. \rm Connect nodes with their counterparts in the duplicate.}\label{2split}
\end{figure}

\begin{figure}
\includegraphics[width=.9\textwidth]{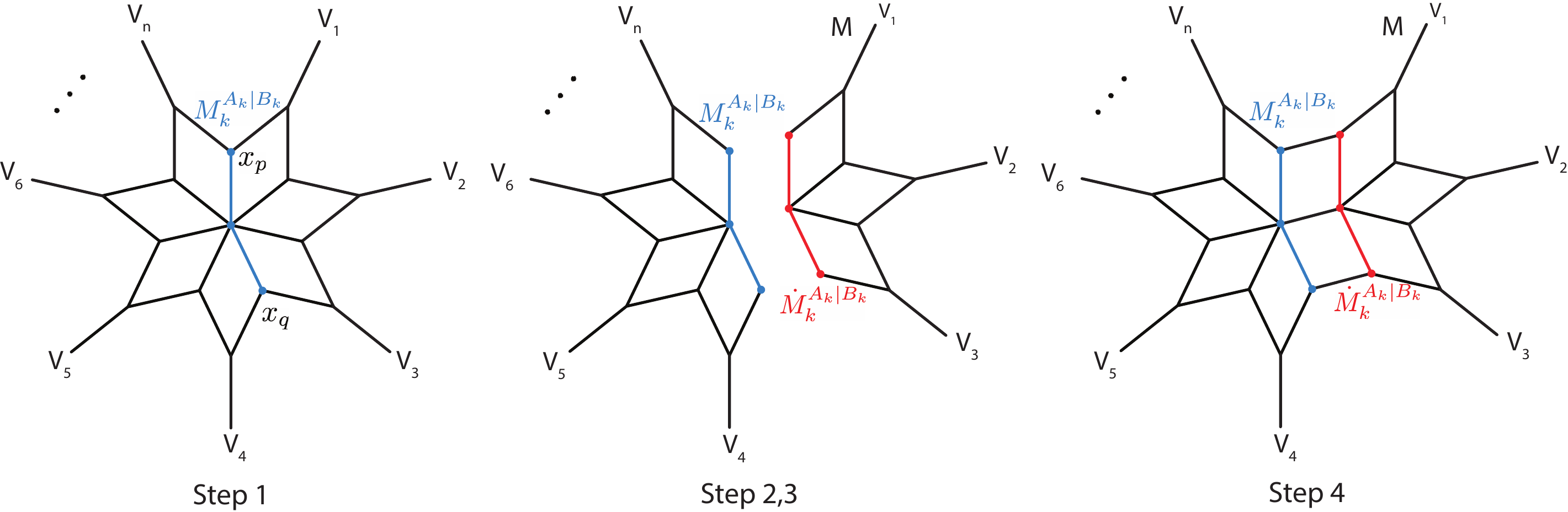}
\caption{Adding the split $S_n^{\{V_1, V_2, V_3\}|\{V_4, \ldots, V_n\}}$ for the set of taxa $\mathcal{X} = \{V_1, \ldots, V_n\}$.
\it Step 1. \rm Determine the split path $M_k^{A_k|B_k}$.
\it Step 2. \rm Create a copy $\dot M_k^{A_k|B_k}$ of $M_k^{A_k|B_k}$.
\it Step 3. \rm Redirect nodes.
\it Step 4. \rm Connect nodes with their counterparts in the duplicate.}\label{2splitB}
\end{figure}

\begin{figure}
\includegraphics[width=\textwidth]{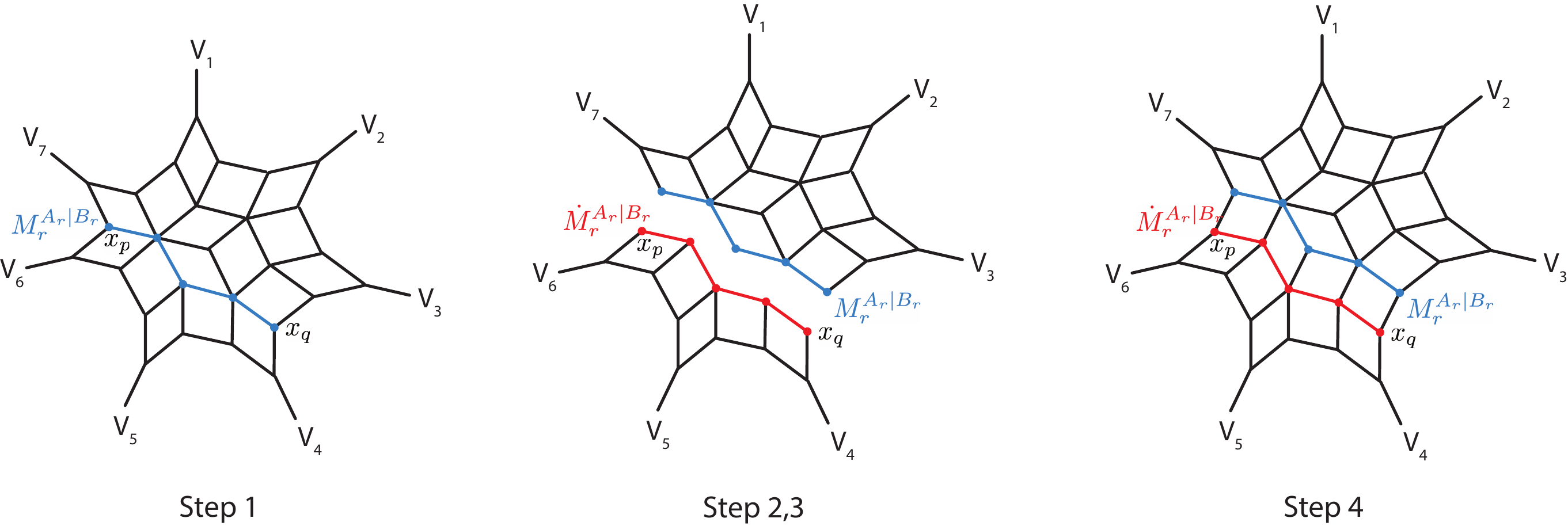}
\caption{Adding the split $S_7^{\{V_4, V_5, V_6\}|\{V_1, V_2, V_3, V_7\}}$ for the set of taxa $\mathcal{X} = \{V_1, \ldots, V_7\}$.
\it Step 1. \rm Determine the split path $M_r^{A_r|B_r}$.
\it Step 2. \rm Create a copy $\dot M_r^{A_r|B_r}$ of $M_r^{A_r|B_r}$.
\it Step 3. \rm Redirect nodes.
\it Step 4. \rm Connect nodes with their counterparts in the duplicate.}\label{2splitC}
\end{figure}

In this section, we introduce the split network constructing algorithm, which produces a planar graph representing the weighted split information contained in a phylogenetic network. Note that the structure of the resultant network depends on the order of the splits being added \cite{huson2010phylogenetic}. We first present the definition of a split path, which is used by the split network constructing algorithm to realize phylogenetic networks. The algorithm takes a set of splits and generates the corresponding phylogenetic network. A few examples are given throughout the rest of the paper: Figure \ref{2split} illustrates how the algorithm adds the first split, Figures \ref{2splitB} and \ref{2splitC} illustrate the process of adding the split $S_i^{A_i|B_i}$ to an existing set of splits $\mathcal{S}^\mathcal{X}$, extending the graph realization of a phylogenetic network with the set of splits $\mathcal{S}^\mathcal{X}$  to a phylogenetic network with the set of splits $\mathcal{S}^\mathcal{X} \cup \left\{S_i^{A_i|B_i}\right\}$. The bottom rows of Figures \ref{fig:N4-stable} and \ref{fig:N5-stable} illustrate the process of extending networks with all possible splits for $n$ = 4 and 5, respectively.

\begin{definition}[Split Path]\label{def:path}\rm
Fix a circular order for a set of taxa of a phylogenetic network $N^\mathcal{X}_x$ with $n$ leaves to be $\mathcal{X} = \{V_1, \dots, V_n\}$ with $V_1 = V_{n+1}$. Let $A_j = \{V_i, \ldots, V_{i+k-1}\}$ and $B_j = \{V_{i+k}, \ldots, V_{i-1}\}$. Let $x_p$ be the node on the shortest path between $V_{i-1}$ and $V_{i}$ on the convex hull of $N^\mathcal{X}_x$, such that the number of edges between $x_p$ and $V_{i-1}$ on the convex hull is the same as the number of edges between $x_p$ and $V_i$ on the convex hull; if such a node does not exist, i.e., the number of edges on the shortest path between $V_{i-1}$ and $V_i$ on the convex hull is odd, then we select $x_p$ to be the node such that the number of edges between $x_p$ and $V_i$ is one edge less than the number of edges between $x_p$ and $V_{i-1}$ on the convex hull. Let $x_q$ be the node on the shortest path between $V_{i+k-1}$ and $V_{i+k}$ on the convex hull of $N^\mathcal{X}_x$, such that the number of edges between $x_q$ and $V_{i+k-1}$ on the convex hull is the same as the number of edges between $x_q$ and $V_{i+k}$ on the convex hull. If such $x_q$ does not exist, then pick $x_q$ to be the node such that the number of edges between $x_q$ and $V_{i+k-1}$ on the convex hull is one edge less than the number of edges between $x_q$ and $V_{i+k}$. The \it split path \rm of the split $S_j^{A_j|B_j}$ is the shortest path between $x_p$ and $x_q$ in $N^\mathcal{X}_x$, denoted by $M_j^{A_j|B_j}$.
\end{definition}

Huson et al.~ introduced the Circular Network Algorithm to generate split networks given a set of weighted splits \cite{huson2010phylogenetic}. In the following, we introduce this algorithm with specified details to fix the number of edges in a phylogenetic network generated by this algorithm. The split network algorithm generates the phylogenetic network $N^{\mathcal{X}}_x(\mathcal{S}^\mathcal{X}_x, \mathcal{W}^\mathcal{X}_x)$ starting from an $n$-star. Then the algorithm iterates elements in $\mathcal{S}^\mathcal{X}_x$, with each iteration adding a split $S_j^{A_j|B_j} \in \mathcal{S}^\mathcal{X}_x$ using the Circular Network Algorithm. Now we present the formal procedure of the split network algorithm as Algorithm \ref{CircularNetworkAlgorithm}.

\begin{algorithm}
\caption{Split Network Algorithm}\label{CircularNetworkAlgorithm}
\begin{algorithmic}[1]
\Procedure{Circular Network Algorithm }{$S_j^{A_j|B_j}$} \Comment{Add the split $S_j^{A_j|B_j}$.}
	\State Determine the split path $M_j^{A_j|B_j}$ leading from $x_p$ to $x_q$.
	\State Create a copy $\dot M_j^{A_j|B_j}$ of $M_j^{A_j|B_j}$.
	\State Redirect every node on the $A_j$ side of  $M_j^{A_j|B_j}$ on the split path to $\dot M_j^{A_j|B_j}$.
	\State Connect each node in $M_j^{A_j|B_j}$ with its counterpart in $\dot M_j^{A_j|B_j}$ with an edge of length $w_j^{A_j|B_j}$.
	\EndProcedure	
	\Procedure{Develop Split Network $N^\mathcal{X}_x$ } {$\mathcal{S}^\mathcal{X}_x$}\Comment{Create $N^{\mathcal{X}}_x(\mathcal{S}^\mathcal{X}_x, \mathcal{W}^\mathcal{X}_x)$.}
    	\For{$j = 1$ \textbf{to} $\left|\mathcal{S}^\mathcal{X}_x\right|$}
		\State Circular Network Algorithm $\left(S_j^{A_j|B_j}\right)$  		 \EndFor
\EndProcedure
\end{algorithmic}
\end{algorithm}

Now we give several examples of phylogenetic networks generated by Algorithm \ref{CircularNetworkAlgorithm}. The first example is starting from an $n$-star and adding the first split $S_j^{A_j | B_j}$ with $A_j = \{V_i, \ldots, V_{i+k-1}\}$ and $B_j = \{V_{i+k}, \ldots, V_{i-1}\}$ for $i = 1$ and $k = 2$, illustrated as Figure \ref{2split}.

\begin{enumerate}
\item[Step 1.] In this example,
$A_j = \{V_1, V_2\}$ and $B_j = \{V_3, \ldots, V_n\}$.
We first find the split path $M_j^{A_j|B_j}:$ to determine $x_p$ and $x_q$; note that the only node between $V_1$ and $V_{n}$ on the convex hull of the $n$-star is the origin, and the only node between $V_2$ and $V_3$ is also the origin. Since both $x_p$ and $x_q$ are located at the origin, $M_j^{A_j|B_j}$ is the trivial path, denoted as a blue dot.

\item[Step 2.] Create a copy $\dot M_j^{A_j|B_j}$ of $M_j^{A_j|B_j}$, which is displayed as the red dot in Figure \ref{2split}.

\item[Step 3.]
For all the edges connected to $M_j^{A_j|B_j}$ on the $A_j$ side, we redirect their end nodes on $M_j^{A_j|B_j}$ to the corresponding nodes on $\dot M_j^{A_j|B_j}$. In this example, the edges connected to $M_j^{A_j|B_j}$ on the $A_j$ side are the nodes of $V_1$ and $V_2$ with valence greater than one.
We redirect these nodes from $M_j^{A_j|B_j}$ to
$\dot M_j^{A_j|B_j}$. Note that the effect is grafting the $V_1$ and $V_2$ leaves to $\dot M_j^{A_j|B_j}$.

\item[Step 4.] Connect pairs of nodes in $M_j^{A_j|B_j}$  and $\dot M_j^{A_j|B_j}$ by a new edge with lengths $w_j^{A_j | B_j}$ as illustrated in Step 4. This edge represents the new split $S_j^{A_j | B_j}$.\end{enumerate}

Now we give an example of extending a phylogenetic network by the split $S_k^{A_k|B_k}$ with $A_k = \{V_1, V_2, V_3\}$ and $B_k = \{V_4, \ldots, V_n\}$. The process is displayed as Figure \ref{2splitB}.

\begin{enumerate}
\item[Step 1.] First we find the path $M_k^{A_k|B_k}$. There is only one node $x_p$ between $V_1$ and $V_{n}$ on the convex hull of the current phylogenetic network, and only one node $x_q$ between $V_3$ and $V_4$  on the convex hull. Thus, the split path $M_k^{A_k|B_k}$ between these two nodes is the blue edge path which travels through the center point in Figure \ref{2splitB}.

\item[Step 2,3.] Create a copy $\dot M_k^{A_k|B_k}$ of $M_k^{A_k|B_k}$, displayed as the red edge path.
For every edge connected to $M_k^{A_k|B_k}$ on the $A_k$ side, we redirect its end node on  $M_k^{A_k|B_k}$  to the corresponding node on $\dot M_k^{A_k|B_k}$.

\item[Step 4.] Connect pairs of nodes in 
 $M_k^{A_k|B_k}$ and $\dot M_k^{A_k|B_k}$ by a set of parallel edges with lengths $w_k^{A_k | B_k}$. This set of parallel edges represents the new split $S_k^{A_k | B_k}$.\end{enumerate}

Now we give the last example of a phylogenetic network extended from an existing one by adding a split $S_r^{A_r|B_r}$ with $A_r = \{V_4, V_5, V_6\}$ and $B_r = \{V_7, V_1, V_2, V_3\}$, displayed as Figure \ref{2splitC}.

\begin{enumerate}
\item[Step 1.]  First we find the path $M_r^{A_r|B_r}$. There are two nodes between $V_3$ and $V_4$ on the convex hull of the existing phylogenetic network, and we select the one closer to $V_4$ to be $x_p$ as specified by Definition \ref{def:path}. There is only one node between $V_6$ and $V_7$ on the convex hull, and that will be $x_q$. The split path $M_r^{A_r|B_r}$ is the blue edge path in Figure \ref{2splitC}.

\item[Step 2,3.]  We create a copy $\dot M_r^{A_r|B_r}$ of $M_r^{A_r|B_r}$, which is the red edge path.
For every edge connected to $M_r^{A_r|B_r}$ on the $A_r$ side, we redirect its end node on $M_r^{A_r|B_r}$ to the corresponding node $\dot M_r^{A_r|B_r}$.

\item[Step 4.]  Connect pairs of nodes in $M_r^{A_r|B_r}$ and $\dot M_r^{A_r|B_r}$ by parallel edges with length $w_r^{A_r | B_r}$, which represent the new split $S_r^{A_r | B_r}$.\end{enumerate}

The structure of phylogenetic networks will allow us to explore their connection with real genus 0 stable curves in $\overline{\mathfrak{M}}_{0,n}(\real)$ and
complex high genus stable curves in $\overline{\mathfrak{M}}_{g,n}(\complex)$. We discuss their connection in the following sections.

\subsection{The Space of Network Topologies and Stable Curves in $\overline{\mathfrak{M}}_{g,n}(\mathbb{C})$}\label{sec:network-curve}

To analyze the connection between phylogenetic networks and stable curves, we first investigate the dual graph that is associated with stable curves for arbitrary genus. Recall that we denoted by $C$ a nodal curve, and $D$ a finite set of smooth points of $C$. The normalization of each irreducible component separates the nodes, and the normalization of $C$ is the disjoint union of the normalization of each irreducible component. Let $(C;D)$ be a connected nodal curve with $n$ marked points. We say that $(C;D)$ is stable if it has a finite automorphism group \cite{arbarello2011geometry}.

\begin{figure}[!htb]\centering
 \includegraphics[width=.9\textwidth]{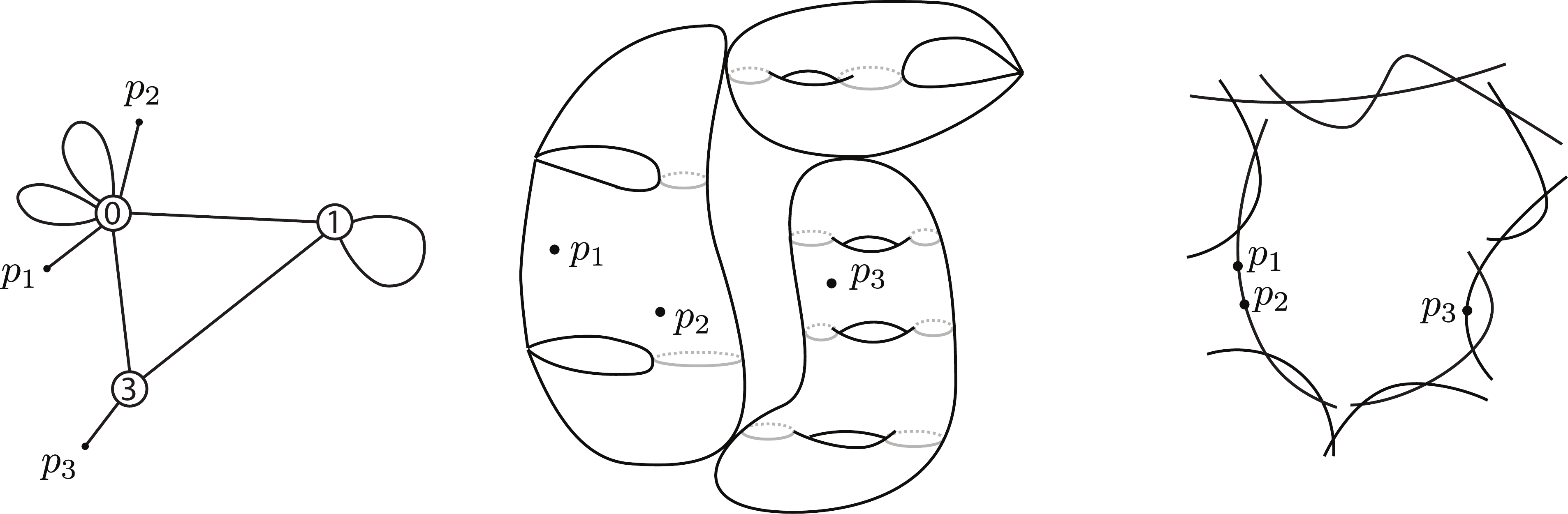}
\caption{A stable curve with 3 marked points and its representations. \it Left: \rm the associated dual graph of the stable curve. \it Middle: \rm a stable curve with 3 marked points and three components. \it Right: \rm The boundary stratum of the stable curve contracted along the gray circles illustrated in the middle graph.}\label{fig:nodal}
\end{figure}

To associate a dual graph Graph$(C;D)$ with a nodal curve $C$ with marked points $D$, we assign a vertex for each component of the normalization of $C$, with the weight of the vertex being the genus of the component. 
The half-edges issuing from a vertex either map to a node of $C$ or to a marked point on the component corresponding to the vertex that this edge issued from. In other words, the edges of the graph are the pairs of half-edges mapping to the same node of $C$, and the legs are the half-edges coming from the marked points. 
This way to associate graphs with stable curves is introduced by Arbarello, Cornalba 􏱥and Griffiths \cite{arbarello2011geometry}, and is referred to as the \it canonical association \rm for the rest of the paper.
Figure \ref{fig:nodal} illustrates a 3-pointed curve (middle) and its dual graph (left). We can pinch the curve along the grey circles in the middle drawing so that it consists of only genus 0 components. The curve contracted along these circles is illustrated in the right-hand drawing, with each piece signifying a copy of $\cp^1$. The contracted curve is on the boundary of the moduli space of the original curve. In particular, the way to associate phylogenetic trees with stable curves is the genus 0 special case of the way to associate phylogenetic networks with stable curves. This is because phylogenetic trees can be viewed as a degeneration of phylogenetic networks by collapsing the sets of parallel edges that induce circles, since the set of pairwise compatible splits is a subset of the set of splits that are compatible with certain circular ordering. Note that although the association between the dual graph of a stable curve and the stable curve is canonical, there is no canonical way to contract a stable curve with high genus components to a curve which has only genus 0 components. Thus, the contraction illustrated in Figure \ref{fig:nodal} from the middle curve with high genus components to the curve on the right that has only genus 0 components is not unique.

\begin{figure}[!htb]\centering
 \includegraphics[width=.8\textwidth]{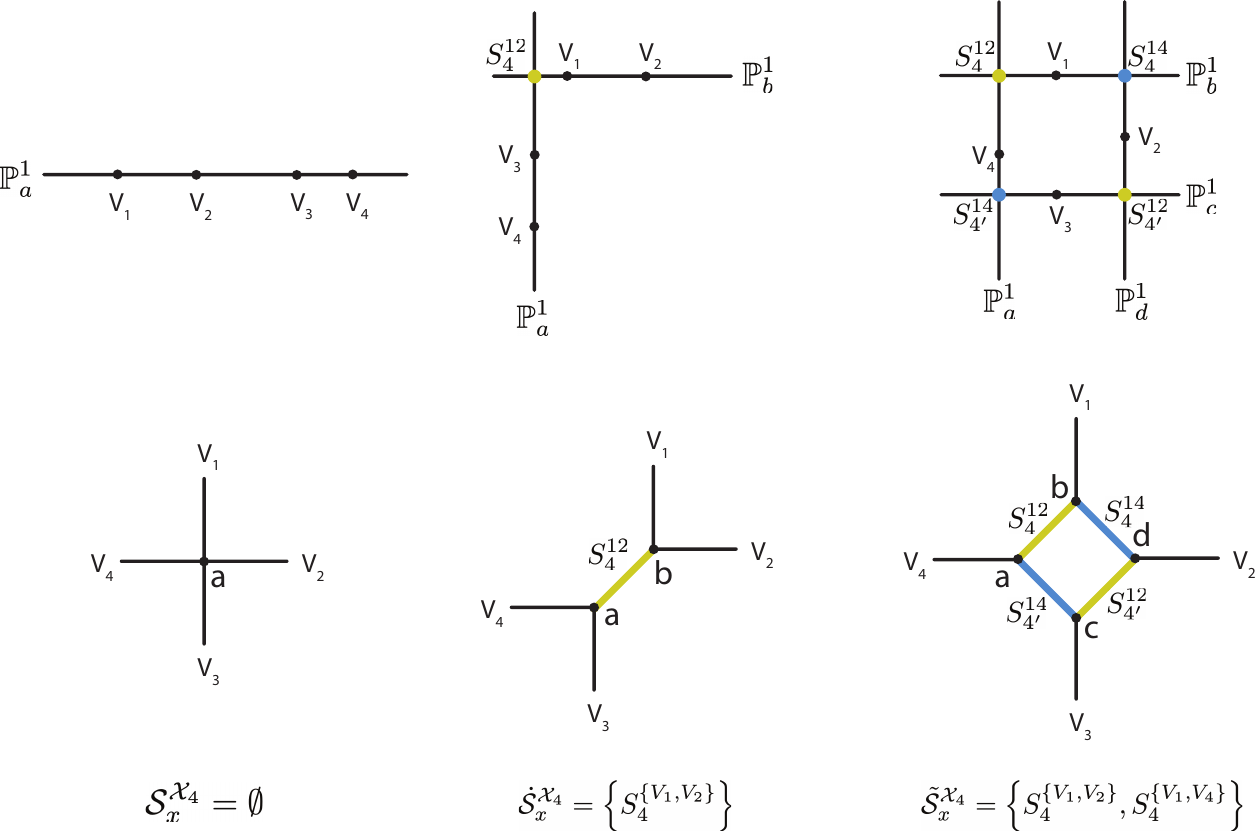}
\caption{Phylogenetic networks with 4 leaves (bottom row) and their counterparts as stable curves in $\overline{\mathfrak{M}}_{0,4}(\mathbb{C})$ and $\overline{\mathfrak{M}}_{1,4}(\mathbb{C})$ (top row). $S_4^{ij}$ represents $S_4^{A_k|B_k}$ with $A_k = \{V_i, V_j\}$ and the concatenation of indices of elements in $A_i$ is smaller than that of $B_i$. \it Left: \rm a boundary point (the origin) $N^{\mathcal{X}_4}_x\left(\emptyset,\emptyset\right)\in\csn_4$ and associated stable curve, an interior point of $\mathfrak{M}_{0,4}(\mathbb{C})$; \it middle: \rm the phylogenetic network $\dot N^{\mathcal{X}_4}_x \left(\dot{\mathcal{S}}^{\mathcal{X}_4}_x,  \dot{\mathcal{W}}^{\mathcal{X}_4}_x  \right)\in\csn_4$ and associated stable curve in $\overline{\mathfrak{M}}_{0,4}(\mathbb{C})$; \it right\rm:  the phylogenetic network $\tilde{N}^{\mathcal{X}_4}_x \left(\tilde{\mathcal{S}}^{\mathcal{X}_4}_x, \tilde{\mathcal{W}}^{\mathcal{X}_4}_x  \right)\in\csn_4$  and associated stable curve on $\partial\overline{\mathfrak{M}}_{1,4}(\mathbb{C})$.
}\label{fig:N4-stable}
\end{figure}

The space of network topologies forms a category and maps injectively into the partition space of high genus stable curves. Before we proceed to formal proof, we first provide a few examples. Figure \ref{fig:N4-stable} displays a phylogenetic network with four vertices and their counterparts as high-genus stable curves by associating each node with a genus 0 component and each edge with an intersection of two components. 
Since we associate each node of the graph of a phylogenetic network with a genus 0 component with no contraction involved, this association between graphs and stable curves is canonical  \cite{arbarello2011geometry}, as introduced at the beginning of Section \ref{sec:network-curve}. 

The phylogenetic networks displayed in the bottom row of Figure \ref{fig:N4-stable} are generated by Algorithm \ref{CircularNetworkAlgorithm}, which adds a new split at each iteration, and their associated stable curves are displayed in the upper row.
The upper-left picture depicts an interior point of $\mathfrak{M}_{0,4}(\mathbb{C})$ with four marked points $V_1$, $V_2$, $V_3$, and $V_4$; it corresponds to a degenerated phylogenetic network $N^{\mathcal{X}_4}_x\left(\mathcal{S}_x^{\mathcal{X}_4},\mathcal{W}_x^{\mathcal{X}_4}\right)\in\csn_4$  as a $4$-star with $\mathcal{S}_x^{\mathcal{X}_4} = \emptyset$, $\mathcal{W}_x^{\mathcal{X}_4} = \emptyset$, and $\mathcal{X}_4 =\{V_1, V_2, V_3, V_4\}$.
 Colliding $V_1$ and $V_2$ sprouts another $\mathbb{P}^1$ with these two points and a special point at the intersection, displayed as the upper-middle picture. This stable curve is on $\partial \overline{\mathfrak{M}}_{0,4}$, which corresponds to a phylogenetic network $\dot N^{\mathcal{X}_4}_x \left(\dot{\mathcal{S}}^{\mathcal{X}_4}_x,  \dot{\mathcal{W}}^{\mathcal{X}_4}_x  \right)\in\csn_4$ with $\dot{\mathcal{S}}^{\mathcal{X}_4}_x = \left\{S_4^{\{V_1, V_2\}|\{V_3, V_4\}}\right\}$
 obtained from adding the split $S_4^{\{V_1, V_2\}|\{V_3, V_4\}}$ into the set of splits $\mathcal{S}_x^{\mathcal{X}_4}$.
 To create a split between $\{V_1, V_2\}$ and $\{V_3, V_4\}$, we first need to find the split path $M_4^{\{V_1, V_2\}|\{V_3, V_4\}}$, which is the node $a$ on the left of the bottom row of Figure \ref{fig:N4-stable}.
According to the Circular Network Algorithm, we duplicate the split path $M_4^{\{V_1, V_2\}|\{V_3, V_4\}}$ and attach $V_1$ and $V_2$ to the duplicated copy. Then we connect the nodes on the split path with their counterparts in the duplicate by a set of parallel edges. This procedure yields the phylogenetic network in the middle, which is a special case of a phylogenetic network that degenerated into a tree. We can continue to add the split $S_4^{\{V_1, V_4\}|\{V_2, V_3\}}$ to obtain the phylogenetic network $\tilde{N}^{\mathcal{X}_4}_x \left(\tilde{\mathcal{S}}^{\mathcal{X}_4}_x, \tilde{\mathcal{W}}^{\mathcal{X}_4}_x  \right)\in\csn_4$ with $\tilde{\mathcal{S}}^{\mathcal{X}_4}_x = \left\{S_4^{\{V_1, V_2\}|\{V_3, V_4\}},S_4^{\{V_1, V_4\}|\{V_2, V_3\}}\right\}$ obtained by adding the split $S_4^{\{V_1, V_4\}|\{V_2, V_3\}}$ into the set of splits $\dot{\mathcal{S}}_x^{\mathcal{X}_4}$. To create a split between $\{V_1, V_4\}$ and $\{V_2, V_3\}$, we first need to find the split path $M_4^{\{V_1, V_4\}|\{V_2, V_3\}}$, which is  the yellow edge between nodes $a$ and $b$ in the middle of the bottom row of Figure \ref{fig:N4-stable}.
Using the Circular Network Algorithm, we duplicate the the split path $M_4^{\{V_1, V_4\}|\{V_2, V_3\}}$ and attach $V_1$ and $V_4$ to the duplicated copy. Then we connect the nodes on the the split path with their counterparts in the duplicate by a set of parallel edges. This procedure yields the phylogenetic network on the right. This induces a circle in the phylogenetic network that corresponds to a genus 1 stable curve in $\partial\overline{\mathfrak{M}}_{1,4}(\complex)$ by associating the nodes $a$, $b$, $c$, and $d$ of the phylogenetic network with genus 0 components $\mathbb{P}^1_a$, $\mathbb{P}^1_b$, $\mathbb{P}^1_c$, and $\mathbb{P}^1_d$ of the stable curve, associating taxa $V_1$, $V_2$, $V_3$, and $V_4$ of the phylogenetic network with marked points $V_1$, $V_2$, $V_3$, and $V_4$ on the stable curve, and associating edges of the phylogenetic network $S_4^{12}$, $S_4^{14}$, $S_{4^\prime}^{12}$, and $S_{4^\prime}^{14}$ with intersections of the zero components of the stable curve.

\begin{figure}[!htb]\centering
 \includegraphics[width=\textwidth]{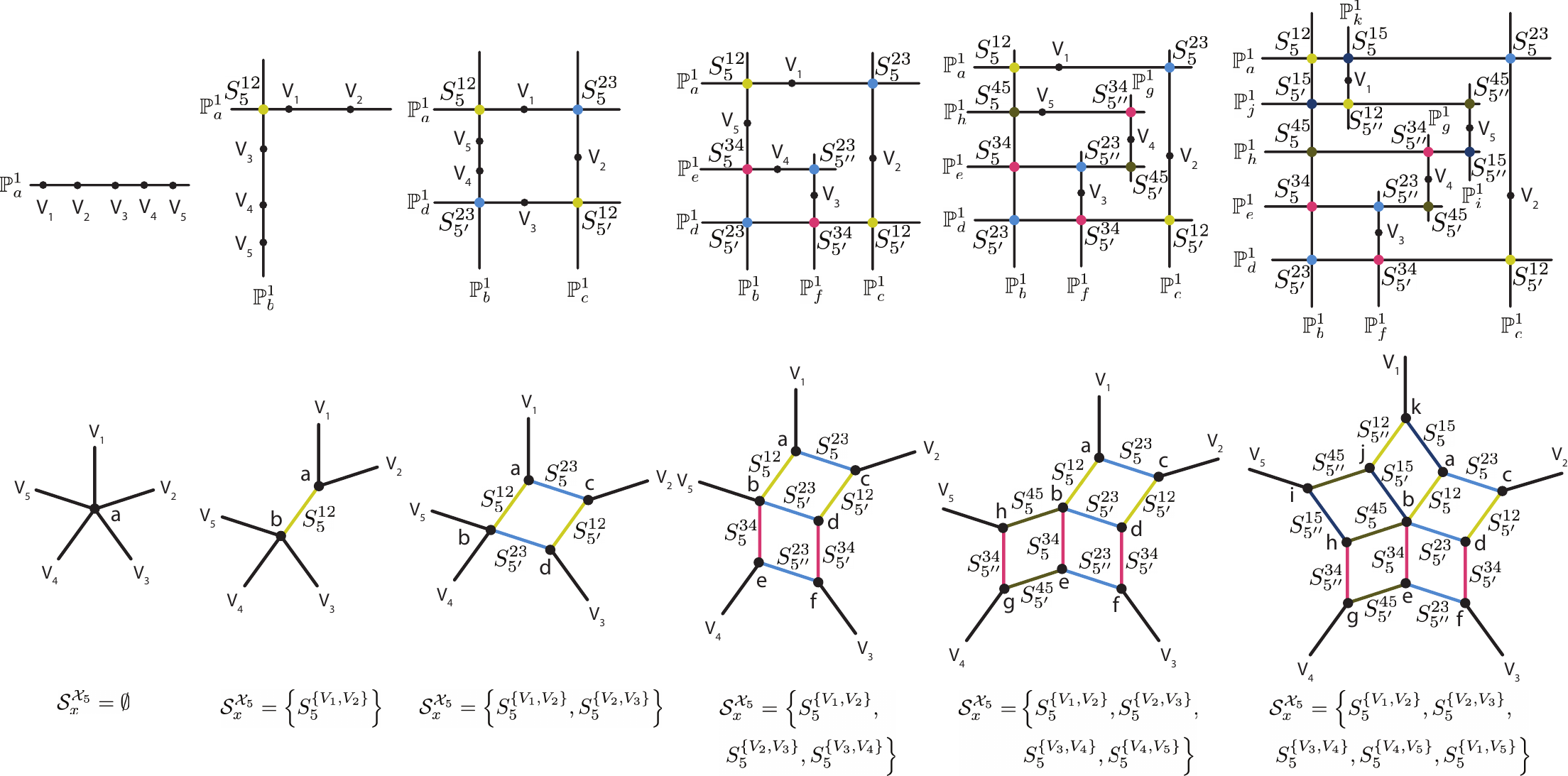}
\caption{
Phylogenetic network with 5 vertices (bottom row) and their counterparts as stable curves in $\overline{\mathfrak{M}}_{g,5}(\mathbb{C})$ for $g = 0,1,2,3,5$ (top row). $S_5^{ij}$ represents $S_5^{A_k|B_k}$ with $A_k = \{V_i, V_j\}$  and the concatenation of indices of elements in $A_i$ is smaller than that of $B_i$.
}\label{fig:N5-stable}
\end{figure}

The case with 5 taxa is more complicated. The number of different splits of a phylogenetic network with $n$ leaves is ${n(n-3)\over 2}$, while the number of edges is a polynomial of the number of leaves (taxa) to the $4^{\text{th}}$-order, i.e., $\mathcal{O}(n^4)$. A detailed discussion on the exact number of edges of a phylogenetic network as a function of the number of leaves is presented in subsequent work.
As illustrated in Figure \ref{fig:N5-stable}, the phylogenetic network in the lower-left contains no split. The next network in the lower-middle contains one split $S_5^{A_1|B_1}$ with $A_1 = \{V_1, V_2\}$ and $B_1 = \{V_3, V_4, V_5\}$. The third phylogenetic network contains two splits, $S_5^{A_1|B_1}$ and $S_5^{A_2|B_2}$ with $A_2 = \{V_2, V_3\}$ and $B_2 = \{V_1, V_4, V_5\}$, and so on.

Stable curves correspond to phylogenetic networks in $\csn_n$ on the boundary of $\overline{\mathfrak{M}}_{g,n}(\mathbb{C})$, obtained from contracting high genus elements into a wedge sum of genus 0 components $\mathbb{P}^1$. In the following, we show that through the phylogenetic network construction algorithm, the graph realization of phylogenetic networks always corresponds to a stable curve. 

\begin{theorem}\rm\label{thm:stable_c}
Phylogenetic networks in $\csn_n$ always correspond to a stable curve in $\overline{\mathfrak{M}}_{g,n}$ by associating each leaf in a phylogenetic network with a marked point, each node that is not on the end of a leaf with a genus 0 component, and each edge as an intersection of two components.
\end{theorem}
\begin{proof}
We proceed by induction on the number of splits for a phylogenetic network with $n$ taxa. The base case is a phylogenetic network $N^\mathcal{X}(\mathcal{S}^\mathcal{X}, \mathcal{W}^\mathcal{X})$ with no internal edge ($\mathcal{S}^\mathcal{X}=\emptyset$); this is an $n$-star. $N^\mathcal{X}$ corresponds to a smooth algebraic curve $\mathcal{C}^X$ with one genus 0 component and a set $X$ of $n$ marked points on the component; $\mathcal{C}^X$ is an interior point in $\mathfrak{M}_{0,n} \subset \overline{\mathfrak{M}}_{\dot g,n}$. For the inductive hypothesis, we assume that a phylogenetic network $\dot N^\mathcal{X}(\dot{\mathcal{S}}^\mathcal{X}, \dot{\mathcal{W}}^\mathcal{X})$ with $|\mathcal{S}^\mathcal{X}| = k$ is associated with a stable curve $\dot{\mathcal{C}}^X$, and we need to prove this is true for $k+1$. Therefore, we show that extending the set of splits $\dot{\mathcal{S}}^\mathcal{X}$ by an arbitrary split $S_n^{A_{k+1}|B_{k+1}}$ that is compatible with the existing circular order using the Split Network Algorithm, the phylogenetic network $\dot N^\mathcal{X}(\dot{\mathcal{S}}^\mathcal{X}, \dot{\mathcal{W}}^\mathcal{X})$ is augmented to a phylogenetic network $\tilde{N}^\mathcal{X}\left(\tilde{\mathcal{S}}^\mathcal{X}, \tilde{\mathcal{W}}^\mathcal{X}\right)$ with $\tilde{\mathcal{S}}^\mathcal{X} = \dot{\mathcal{S}}^\mathcal{X} \cup \left\{S_n^{A_{k+1}|B_{k+1}}\right\}$, such that $\tilde{N}^\mathcal{X}$ is still associated with a stable curve $\tilde{\mathcal{C}}^X$ in $\overline{\mathfrak{M}}_{\tilde g,n}$.

To show that the algebraic curve $\tilde{\mathcal{C}}^X$ is stable, we need to investigate whether it has finite automorphism groups on each component, i.e., whether each component has at least three special points. 
Special points on an algebraic curve are the marked points and intersections which correspond to leaves and edges of phylogenetic networks, respectively; components of an algebraic curve correspond to internal nodes of  the phylogenetic network. Therefore, whether an algebraic curve is stable can be determined by whether each internal node of the phylogenetic network associated with this algebraic curve has at least valence three. The inductive hypothesis assumes that each component of the stable curve $\dot{\mathcal{C}}^{X}$ has at least three special points. Using Algorithm \ref{CircularNetworkAlgorithm}, adding a split $S_n^{A_{k+1}|B_{k+1}}$ only affects nodes along the split path $M_n^{A_{k+1}|B_{k+1}}$ separating $A_{k+1}$ and $B_{k+1}$, and the connectivity of nodes not on $M_n^{A_{k+1}|B_{k+1}}$ remains the same. Therefore, we only need to examine whether the components of $\tilde{\mathcal{C}}^X$ corresponding to the nodes along  $M_n^{A_{k+1}|B_{k+1}}$ and its duplicate $\dot M_n^{A_{k+1}|B_{k+1}}$ have at least three special points.

However, each node on either copy of the split path has at least two edges connected to this node, which are the edges connected to this node on the path. Then we connect nodes on  $M_n^{A_{k+1}|B_{k+1}}$ with their counterparts on the duplicated path $\dot M_n^{A_{k+1}|B_{k+1}}$ with a set of parallel edges. This gives a third edge connected to each node on the split path as well as on its duplicate. Thus, all nodes in $\tilde{N}^\mathcal{X}(\tilde{\mathcal{S}}^\mathcal{X}, \tilde{\mathcal{W}}^\mathcal{X})$ have at least valence three, and we conclude that the $\tilde{\mathcal{C}}^X$ is stable and is contained in $\overline{\mathfrak{M}}_{\tilde g,n}$.

\end{proof}

We define the category of network $\mathfrak{N}_n$ with objects $\mathcal{O}(\mathfrak{N}_n)$ being the collection of networks $N^{\mathcal{X}}_x(\mathcal{S}^\mathcal{X}_x, \mathcal{W}^\mathcal{X}_x)$ in $\csn_n$ with the set of splits $\mathcal{S}^\mathcal{X}_x$ and the set of internal edge length $\mathcal{W}^\mathcal{X}_x$.
For two objects $N^\mathcal{X}_x(\mathcal{S}^\mathcal{X}_x, \mathcal{W}^\mathcal{X}_x),N^\mathcal{X}_y(\mathcal{S}^\mathcal{X}_y, \mathcal{W}^\mathcal{X}_y) \in \mathcal{O}(\mathfrak{N}_n)$, we define a morphism $\mor(N^{\mathcal{X}}_x, N^{\mathcal{X}}_y)$ between $N^{\mathcal{X}}_x$ and $N^{\mathcal{X}}_y$ if $\mathcal{S}^\mathcal{X}_y \subset \mathcal{S}^\mathcal{X}_x$. If $\mathcal{S}^\mathcal{X}_x = \mathcal{S}^\mathcal{X}_y$, we define the morphism $\mor(N^{\mathcal{X}}_x, N^{\mathcal{X}}_y)$ to be the identity.

\begin{proposition}\rm\label{thm:cat_network}
$\mathfrak{N}_n$ is a category.
\end{proposition}

\begin{proof}
By definition, $\mor(N^{\mathcal{X}}_x, N^{\mathcal{X}}_x)$ is the identity. Therefore, we have an identity morphism $\id_{N^{\mathcal{X}}_x}\in \mathcal{M}(\mathfrak{N}_n)$ for each object $N^{\mathcal{X}}_x$. We need to show that we have the composition law:
$$\circ: \mor( N^{\mathcal{X}}_x, N^{\mathcal{X}}_y ) \times \mor( N^{\mathcal{X}}_y, N^{\mathcal{X}}_z )\longrightarrow\mor( N^{\mathcal{X}}_x , N^{\mathcal{X}}_z )$$
for each triple of objects $N^{\mathcal{X}}_x$, $N^{\mathcal{X}}_y$, $N^{\mathcal{X}}_z$. This is true because if $\mathcal{S}^\mathcal{X}_y \subset \mathcal{S}^\mathcal{X}_x$ and $\mathcal{S}^\mathcal{X}_z \subset \mathcal{S}^\mathcal{X}_y$, then $\mathcal{S}^\mathcal{X}_z \subset \mathcal{S}^\mathcal{X}_x$. The associativity of composition is true due to the set relation:
$$
(\mor ( N^{\mathcal{X}}_x, N^{\mathcal{X}}_y ) \circ \mor ( N^{\mathcal{X}}_y, N^{\mathcal{X}}_z ))  \circ \mor ( N^{\mathcal{X}}_z, N^{\mathcal{X}}_w ) \longrightarrow
\mor ( N^{\mathcal{X}}_x, N^{\mathcal{X}}_z )
 \circ \mor ( N^{\mathcal{X}}_z, N^{\mathcal{X}}_w ) 
 \longrightarrow
\mor ( N^{\mathcal{X}}_x, N^{\mathcal{X}}_w),
$$
and also
$$
\mor ( N^{\mathcal{X}}_x, N^{\mathcal{X}}_y ) \circ (\mor ( N^{\mathcal{X}}_y, N^{\mathcal{X}}_z ) \circ \mor ( N^{\mathcal{X}}_z, N^{\mathcal{X}}_w ) )
 \longrightarrow
\mor ( N^{\mathcal{X}}_x, N^{\mathcal{X}}_y ) \circ 
\mor ( N^{\mathcal{X}}_y, N^{\mathcal{X}}_w)
 \longrightarrow
\mor ( N^{\mathcal{X}}_x, N^{\mathcal{X}}_w).
$$
The identity
$$\mor ( N^{\mathcal{X}}_x, N^{\mathcal{X}}_x ) \circ \mor ( N^{\mathcal{X}}_x, N^{\mathcal{X}}_y ) =\mor ( N^{\mathcal{X}}_x, N^{\mathcal{X}}_y ),$$
and
$$\mor ( N^{\mathcal{X}}_x, N^{\mathcal{X}}_y ) \circ \mor ( N^{\mathcal{X}}_x, N^{\mathcal{X}}_x ) =\mor ( N^{\mathcal{X}}_x, N^{\mathcal{X}}_y )$$
also follows.
\end{proof}

\subsection{Categorical Equivalence of Phylogenetic Networks and Boundary Divisors in $\overline{\mathfrak{M}}_{0,n}(\mathbb{R})$}\label{subseq:eq_real}
In the previous section, we have defined the category of phylogenetic networks and its inclusion relation into the boundary of $\overline{\mathfrak{M}}_{g,n}(\mathbb{C})$. In this section, we present the categorical equivalence between the category of phylogenetic networks and the category of boundary divisors of $\overline{\mathfrak{M}}_{0,n}(\mathbb{R})$. 
We first define the category of boundary divisors of $\overline{\mathfrak{M}}_{0,n}(\mathbb{R})$.
Recall that a \it divisor \rm is an element of the free abelian group generated by the subvarieties of codimension one \cite{hartshorne2013algebraic}. Restricting to $\mathbb{R}$, Definition \ref{def:boundary} gives that the \it boundary \rm $\Delta = \overline{\mathfrak{M}}_{g,n}(\mathbb{R}) - \mathfrak{M}_{g,n}(\mathbb{R})$ of $\overline{\mathfrak{M}}_{g,n}(\mathbb{R}) $ is a divisor, with each component the closure of a locus of curves with 1 node. Furthermore, according to Definition \ref{fulton1996notes}, the \it boundary divisor \rm of $\overline{\mathfrak{M}}_{0,n}(\mathbb{R})$ corresponding to the marking partition $A \cup B = [n]$ is isomorphic to
$$\overline{\mathfrak{M}}_{0,A\cup\{\cdot\}}(\mathbb{R}) \times \overline{\mathfrak{M}}_{0,B\cup\{\cdot\}}(\mathbb{R}).$$

The marked points $A_i$ and $B_i$ on each component of the boundary divisors of $\overline{\mathfrak{M}}_{g,n}$ induce a bipartition of the set of marked points $X$ with $|X|=n$ such that $A_i \cup B_i = X$, $A_i \cap B_i = \emptyset$, $|A_i|\geq 2$, and $|B_i|\geq 2$ by the stable condition.
We define the category of boundary divisors of $\overline{\mathfrak{M}}_{0,n}(\mathbb{R})$ with objects $\mathcal{O}(\mathfrak{B}_n)$ being sets of boundary divisors of $\overline{\mathfrak{M}}_{0,n}(\mathbb{R})$, where each set contains boundary divisors that induce bipartitions that are compatible with a certain circular ordering. Given an element $\mathcal{B}^X_x \in \mathcal{O}(\mathfrak{B}_n)$, the compatibility criteria is equivalent to any two boundary divisors in $\mathcal{B}^X_x$ that induce partition $A_i \cup B_i$ and $A_j \cup B_j$, which satisfies:
$$D_n(A_i \sqcup B_i) = D_{|A_j|}(A_j) \sqcup D_{|B_j|}(B_j),$$
where $D_k \in \mathrm{Dih}_{2k}$ and $\sqcup$ represents concatenation. We give an example of two boundary divisors that share the same circular ordering as the following:\\

\noindent \bf Example. \rm Set $A_i = (V_2, V_3, V_4)$, $B_i = (V_5, V_1)$, $A_j = (V_3, V_4)$, and $B_j = (V_2, V_1, V_5)$. Then 
$$D_5(A_i \sqcup B_i) = D_5(V_2, V_3, V_4, V_5, V_1) = (V_3, V_4, V_5, V_1, V_2)$$
for a rotational element $D_5$ in the Dihedral group $\mathrm{Dih}_{10}$. We also have
$$D_2(A_j) \sqcup D_3(B_j) = D_2(V_3, V_4) \sqcup D_3(V_2, V_1, V_5) = (V_3, V_4, V_5, V_1, V_2)$$
for $D_2$ being the identity element in $\mathrm{Dih}_4$ and $D_3$ being the reflection element in $\mathrm{Dih}_6$. This is an example of two boundary divisors compatible with the circular order $(V_3, V_4, V_5, V_1, V_2)$. By definition, all the boundary divisors in an element of $\mathcal{O}(\mathfrak{B}_n)$ are compatible with a fixed circular ordering.\\

We denote the set of marking partitions $A_i \cup B_i = X$ induced by the set of boundary divisors of an element $\mathcal{B}^X_x \in \mathcal{O}(\mathfrak{B}_n)$ as
$$\mathcal{S}^X_x = \left\{S^{A_1\mid B_1}_{n}, \ldots, S^{A_k\mid B_k}_{n}\right\}.$$
We define a morphism $\mor(\mathcal{B}^X_x, \mathcal{B}^X_y) \in\mathcal{M}(\mathfrak{B}_n)$ between two elements $\mathcal{B}^X_x, \mathcal{B}^X_y \in \mathcal{O}(\mathfrak{B}_n)$ if 
 the set of marking partitions $\mathcal{S}^\mathcal{X}_y \subset \mathcal{S}^\mathcal{X}_x$.

\begin{proposition}\rm\label{thm:cat_curve}
$\mathfrak{B}_n$ is a category.
\end{proposition}
\begin{proof}
To show that $\mathfrak{B}_n$ is a category, first we note that we have an identity morphism $\id_{\mathcal{B}^X_x} = \mor(\mathcal{B}^X_x,\mathcal{B}^X_x) \in\mathcal{M}(\mathfrak{B}_n)$ for each object $\mathcal{B}^X_x$ since $\mathcal{S}^X_x \subset \mathcal{S}^X_x$. We also have the composition law because for each triple of objects $\mathcal{B}^X_x$, $\mathcal{B}^X_y$, $\mathcal{B}^X_z$, such that $\mathcal{S}^X_y \subset \mathcal{S}^X_x$ and $\mathcal{S}^X_z \subset \mathcal{S}^X_y$, we have $\mathcal{S}^X_z \subset \mathcal{S}^X_x$. Thus, $$\circ: \mor ( \mathcal{B}^X_x , \mathcal{B}^X_y ) \times \mor ( \mathcal{B}^X_y , \mathcal{B}^X_z )\longrightarrow\mor(\mathcal{B}^X_x ,\mathcal{B}^X_z ).$$

Given $\mathcal{S}^X_y \subset \mathcal{S}^X_x$,  $\mathcal{S}^X_z \subset \mathcal{S}^X_y$, and $\mathcal{S}^X_w \subset \mathcal{S}^X_z$, the composition is also associative: $\mathcal{S}^X_y \subset \mathcal{S}^X_x$,  $\mathcal{S}^X_z \subset \mathcal{S}^X_y$ gives $\mathcal{S}^X_z \subset \mathcal{S}^X_x$. This combined with $\mathcal{S}^X_w \subset \mathcal{S}^X_z$ gives $\mathcal{S}^X_w \subset \mathcal{S}^X_x$. Therefore, we have:
$$(\mor ( \mathcal{B}^X_x , \mathcal{B}^X_y ) \times \mor ( \mathcal{B}^X_y , \mathcal{B}^X_z )) \times \mor(\mathcal{B}^X_z ,\mathcal{B}^X_w ) \longrightarrow \mor(\mathcal{B}^X_x ,\mathcal{B}^X_w ).$$
This is equal to the composition switching the order of the association, since $\mathcal{S}^\mathcal{X}_z \subset \mathcal{S}^\mathcal{X}_y$ and  $\mathcal{S}^\mathcal{X}_w \subset \mathcal{S}^\mathcal{X}_z$ imply $\mathcal{S}^\mathcal{X}_w \subset \mathcal{S}^\mathcal{X}_y$; combining with $\mathcal{S}^\mathcal{X}_y \subset \mathcal{S}^\mathcal{X}_x$, we have $\mathcal{S}^\mathcal{X}_w \subset \mathcal{S}^\mathcal{X}_x$.
$$ \mor ( \mathcal{B}^X_x , \mathcal{B}^X_y ) \times (\mor ( \mathcal{B}^X_y , \mathcal{B}^X_z ) \times \mor(\mathcal{B}^X_z ,\mathcal{B}^X_w) )
\longrightarrow \mor ( \mathcal{B}^X_x , \mathcal{B}^X_w ).$$
Finally, we need to check the identity morphisms. For any identity morphism $\mor ( \mathcal{B}^X_x , \mathcal{B}^X_x )$ and any morphism $\mor ( \mathcal{B}^X_y , \mathcal{B}^X_z ) \in \mathcal{M}(\mathfrak{B}_n)$,
$$ \mor ( \mathcal{B}^X_x , \mathcal{B}^X_x ) \circ \mor ( \mathcal{B}^X_y , \mathcal{B}^X_z )\longrightarrow\mor ( \mathcal{B}^X_y , \mathcal{B}^X_z ),$$
because $\mathcal{S}^\mathcal{X}_x \subset \mathcal{S}^\mathcal{X}_x$ and $\mathcal{S}^\mathcal{X}_z \subset \mathcal{S}^\mathcal{X}_y$ implies $\mathcal{S}^\mathcal{X}_z \subset \mathcal{S}^\mathcal{X}_y$. Similarly,
$$\mor ( \mathcal{B}^X_y , \mathcal{B}^X_z )\circ \mor ( \mathcal{B}^X_x , \mathcal{B}^X_x )\longrightarrow\mor ( \mathcal{B}^X_y , \mathcal{B}^X_z ).$$
\end{proof}

Now we start our discussion on the functorial relationship between the category of phylogenetic networks $\mathfrak{N}_n$ and the category of boundary divisors $\mathfrak{B}_n$. We first introduce a forgetful functor $H$ mapping from the category of phylogenetic networks $\mathfrak{N}_n$ to the category of boundary divisors $\mathfrak{B}_n$. The functor $H$ maps onto the category of phylogenetic networks, forgetting the weights of internal edges of each phylogenetic network in $\mathfrak{N}_n$.

\begin{proposition}\rm\label{prop:c_to_t}
The functor 
\begin{align*}
H : \mathfrak{N}_n & \longrightarrow \mathfrak{B}_n \\
N^\mathcal{X}_x(\mathcal{S}^\mathcal{X}_x, \mathcal{W}^\mathcal{X}_x)
&\mapsto \mathcal{B}^X_x(\mathcal{S}^X_x) \\
\mor(N^\mathcal{X}_x,N^\mathcal{X}_y)
&\mapsto \mor\left(\mathcal{B}^X_x,\mathcal{B}^X_y\right)
\end{align*}
maps from the category of networks $\mathfrak{N}_n$
onto the category of boundary divisors $\mathfrak{B}_n$.
\end{proposition}
\begin{proof}
The functor $H$ assigns each object $N^\mathcal{X}_x(\mathcal{S}^\mathcal{X}_x, \mathcal{W}^\mathcal{X}_x)\in \mathcal{O}(\mathfrak{N}_n)$ an object $\mathcal{B}^X_x(\mathcal{S}^X_x)\in \mathcal{O}(\mathfrak{B}_n)$. This assignment is well-defined because elements in a set of bipartitions $\mathcal{S}^\mathcal{X}_x$ of a phylogenetic network $N^\mathcal{X}_x$ share the same circular ordering, so there exists a set of boundary divisors $\mathcal{B}^X_x(\mathcal{S}^X_x)$ which induce the set of bipartitions $\mathcal{S}^X_x$, and $N^\mathcal{X}_x(\mathcal{S}^\mathcal{X}_x, \mathcal{W}^\mathcal{X}_x)$ is mapped to $\mathcal{B}^X_x(\mathcal{S}^X_x)$.
The functor $H$ assigns to each morphism 
$\mor(N^\mathcal{X}_x, N^\mathcal{X}_y) \in \mathcal{M}(\mathfrak{N}_n)$ a morphism
$H(\mor(N^\mathcal{X}_x, N^\mathcal{X}_y) )= \mor\left(\mathcal{B}^X_x,\mathcal{B}^X_y\right)\in \mor\left(\mathfrak{B}_n\right)$, such that
$$H\left(\id_{N^\mathcal{X}_x}\right) =H\left(\mor(N^\mathcal{X}_x, N^\mathcal{X}_x)\right)
=\mor\left(\mathcal{B}^X_x, \mathcal{B}^X_x\right)
 = \id_{H(N^\mathcal{X}_x)} = \id_{\mathcal{B}^X_x}.$$
Now we consider two morphisms
$\mor(N^\mathcal{X}_x,N^\mathcal{X}_y), \mor(N^\mathcal{X}_y,N^\mathcal{X}_z) \in \mathcal{M}(\mathfrak{N}_n)$, such that
$$H\left(\mor(N^\mathcal{X}_x, N^\mathcal{X}_y)\circ \mor(N^\mathcal{X}_y, N^\mathcal{X}_z)\right)
= H\left(\mor(N^\mathcal{X}_x, N^\mathcal{X}_z)\right)=\mor\left(\mathcal{B}^X_x,\mathcal{B}^X_z\right).$$
Also,
$$H\left(\mor(N^\mathcal{X}_x, N^\mathcal{X}_y)\right)\circ H\left(\mor(N^\mathcal{X}_y, N^\mathcal{X}_z)\right)
=\mor\left(\mathcal{B}^X_x,\mathcal{B}^X_y\right) \circ \mor\left(\mathcal{B}^X_y,\mathcal{B}^X_z\right) 
=\mor\left(\mathcal{B}^X_x,\mathcal{B}^X_z\right).$$
Thus, we have
$$H\left(\mor(N^\mathcal{X}_x, N^\mathcal{X}_y)\circ \mor(N^\mathcal{X}_y, N^\mathcal{X}_z)\right)=H\left(\mor(N^\mathcal{X}_x, N^\mathcal{X}_y)\right)\circ H\left(\mor(N^\mathcal{X}_y, N^\mathcal{X}_z)\right),$$
so $H$ is a functor from $\mathfrak{N}_n$ to $\mathfrak{B}_n$. Because elements in the set of bipartitions $\mathcal{S}^X_x$ induced by $\mathcal{B}^X_x\in\mathcal{O}(\mathfrak{B}_n)$ share the same cyclic ordering by definition, they can always find (non-unique) preimages in the category $\mathfrak{N}_n$. Thus, the functor $H$ maps the set of objects in the domain category onto the set of objects in the range category.
\end{proof}

The category of networks and the category of boundary divisors are connected through network topologies. We define the category of network topologies $\overline{\mathfrak{N}}_n$ with objects $\mathcal{O}(\overline{\mathfrak{N}}_n)$ as the collection of equivalence classes of networks $\overline{N}^\mathcal{X}_x(\mathcal{S}^\mathcal{X}_x)$ in $\csn_n$ that share the same set of splits $\mathcal{S}^\mathcal{X}_x$ with different weights $\mathcal{W}^\mathcal{X}_x$. For two objects $\overline{N}^\mathcal{X}_x(\mathcal{S}^\mathcal{X}_x), \overline{N}^\mathcal{X}_y(\mathcal{S}^\mathcal{X}_y) \in \mathcal{O}(\overline{\mathfrak{N}}_n)$, we define a morphism $\mor\left(\overline{N}^\mathcal{X}_x, \overline{N}^\mathcal{X}_y\right)$ if $\mathcal{S}^\mathcal{X}_y\subset \mathcal{S}^\mathcal{X}_x$ as the following: Let $\mathcal{I} = \{i_1, \ldots, i_r\}$ denote the indices of the set $\mathcal{S}^\mathcal{X}_x \setminus \mathcal{S}^\mathcal{X}_y$; we define a projection map $\pi_i$ that forgets the $i$-th component, and denote by $\pi_{\mathcal{I}}$ the composition of projection maps $\pi_{i_1} \circ \cdots \circ \pi_{i_r}$, then
\begin{align*}
\mor\left(\overline{N}^\mathcal{X}_x, \overline{N}^\mathcal{X}_y\right):\mathcal{O}(\mathfrak{N}_n) &\longrightarrow \mathcal{O}(\mathfrak{N}_n)\\
\overline{N}^\mathcal{X}_x(\mathcal{S}^\mathcal{X}_x) & \mapsto \pi_\mathcal{I} \left(\overline{N}^\mathcal{X}_x\left(\mathcal{S}^\mathcal{X}_x\right)\right)=\overline{N}^\mathcal{X}_y(\mathcal{S}^\mathcal{X}_y) .
\end{align*}
If $\mathcal{S}^\mathcal{X}_x = \mathcal{S}^\mathcal{X}_y$, we define the morphism $\mor\left(\overline{N}^\mathcal{X}_x, \overline{N}^\mathcal{X}_y\right)$ as the identity. 
The proof that $\overline{\mathfrak{N}}_n$ is a category follows from the proof of $\mathfrak{N}_n$ being a category. To show that the category of network topology $\overline{\mathfrak{N}}_n$ and the category of boundary divisors $\mathfrak{B}_n$ are equivalent, we first define two functors between these two categories.

The functor 
\begin{align*}
\overline{H} : \overline{\mathfrak{N}}_n & \longrightarrow \mathfrak{B}_n \\
\overline{N}^\mathcal{X}_x(\mathcal{S}^\mathcal{X}_x)&\mapsto \mathcal{B}^X_x(\mathcal{S}^X_x) \\
\mor\left(\overline{N}^\mathcal{X}_x,\overline{N}^\mathcal{X}_y\right)
&\mapsto \mor\left(\mathcal{B}^X_x, \mathcal{B}^X_y\right)
\end{align*}
maps from the category of network topology $\overline{\mathfrak{N}}_n$ to the category of boundary divisors $\mathfrak{B}_n$ in $\overline{\mathcal{M}}_{0,n}(\mathbb{R})$. The functor
\begin{align*}
\overline{K}: \mathfrak{B}_n & \longrightarrow \overline{\mathfrak{N}}_n \\
\mathcal{B}^X_x(\mathcal{S}^X_x)
&\mapsto\overline{N}^\mathcal{X}_x(\mathcal{S}^\mathcal{X}_x) \\
\mor(\mathcal{B}^X_x,\mathcal{B}^X_y)
&\mapsto \mor\left(\overline{N}^\mathcal{X}_x, \overline{N}^\mathcal{X}_y\right)
\end{align*}
maps from the category  of boundary divisors $\mathfrak{B}_n$ in $\overline{\mathcal{M}}_{0,n}(\mathbb{R})$ to the category of network topology.

\begin{lemma}\rm\label{lemma:iso}
The function $\overline{h}$ induced by the functor $\overline{H}$:
\begin{align*}
\overline{h} : \mathcal{O}(\overline{\mathfrak{N}}_n) & \longrightarrow \mathcal{O}(\mathfrak{B}_n)\\
\overline{N}^\mathcal{X}_x(\mathcal{S}^\mathcal{X}_x)&\mapsto \mathcal{B}^X_x(\mathcal{S}^X_x)\end{align*}
is an isomorphism with inverse $\overline{h}^{-1} = \overline{k}$ induced by the functor $\overline{K}$:
\begin{align*}
\overline{k}: \mathcal{O}(\mathfrak{B}_n) & \longrightarrow \mathcal{O}(\overline{\mathfrak{N}}_n) \\
\mathcal{B}^X_x(\mathcal{S}^X_x)&\mapsto 
\overline{N}^\mathcal{X}_x(\mathcal{S}^\mathcal{X}_x).\end{align*}
\end{lemma}
\begin{proof}
The set of bipartitions $\mathcal{S}^X_x$ representing marking partitions $A_i \cup B_i = [n]$ for $n$ marked points of boundary divisors uniquely determines an element $\mathcal{B}^X_x(\mathcal{S}^X_x)\in\mathcal{O}(\mathfrak{B}_n)$. Similarly, the set of splits $\mathcal{S}^\mathcal{X}_x$ also uniquely specifies an element $\overline{N}^\mathcal{X}_x(\mathcal{S}^\mathcal{X}_x)\in\mathcal{O}(\overline{\mathfrak{N}}_n)$. Therefore, $\bar h$ and $\bar k$ are well defined; to show that they are isomorphisms, we only need to show that there is an isomorphism between sets of bipartitions induced by elements in $ \mathcal{O}(\mathfrak{B}_n)$ and sets of splits for elements in $\mathcal{O}(\overline{\mathfrak{N}}_n)$.

Because $\overline{\mathfrak{M}}_{0,n}(\mathbb{R})$ requires $n\geq 4$ as discussed in Section \ref{sec:intro}, boundary divisors as elements of objects in $\mathcal{O}(\mathfrak{B}_n)$ inherit this requirement. Furthermore, each component of a boundary divisor in $\overline{\mathfrak{M}}_{0,n}(\mathbb{R})$ needs at least 2 marked points to satisfy the stable criteria, i.e., for the marking partition $A\cup B = [n]$, we need $|A|\geq 2$ and $|B|\geq 2$. Similarly, elements in $\mathcal{O}(\overline{\mathfrak{N}}_n)$ also require $n\geq 4$ since $|A_i|\geq 2$ and $|B_i|\geq 2$ for all splits $S_n^{A_i|B_i}$ in $\overline{N}^\mathcal{X}_x$. Therefore, the set of bipartitions induced by boundary divisors of $\overline{\mathfrak{M}}_{0,n}(\mathbb{R})$ is isomorphic to the set of splits of network topologies in $\mathcal{O}(\overline{\mathfrak{N}}_n)$.

Now we need to show that a set of bipartitions $S_n^{A_i|B_i}$ forms an element in $ \mathcal{O}(\mathfrak{B}_n)$ if and only if  the corresponding set of splits forms an element in $ \mathcal{O}(\overline{\mathfrak{N}}_n)$. Note that elements $\mathcal{B}^X_x(\mathcal{S}^X_x)\in\mathcal{O}(\mathfrak{B}_n)$ are sets of boundary divisors such that their marking partitions $S_n^{A_i|B_i} \in \mathcal{S}^X_x$ share the same circular ordering; elements $\overline{N}^\mathcal{X}_x(\mathcal{S}^\mathcal{X}_x)\in\mathcal{O}(\overline{\mathfrak{N}}_n)$ are phylogenetic network topologies such that their splits $S_n^{A_i|B_i}\in\mathcal{S}^\mathcal{X}_x$ share the same circular ordering by definition. Since the induced bipartitions of $\mathcal{B}^X_x(\mathcal{S}^X_x)$ and splits of $\overline{N}^\mathcal{X}_x(\mathcal{S}^\mathcal{X}_x)$ are isomorphic, and they form an object in their respective category following the same condition, the sets of partitions that determine objects in $ \mathcal{O}(\mathfrak{B}_n)$ and $\mathcal{O}(\overline{\mathfrak{N}}_n)$ are isomorphic. Therefore, $\overline{h}$ and $\overline{k}$ are isomorphisms between  $ \mathcal{O}(\mathfrak{B}_n)$ and $\mathcal{O}(\overline{\mathfrak{N}}_n)$.
\end{proof}

\begin{theorem}\rm\label{thm:equiv}
The category $\overline{\mathfrak{N}}_n$ and the category $\mathfrak{B}_n$ are equivalent.
\end{theorem}
\begin{proof}
To show that the category $\overline{\mathfrak{N}}_n$ and the category $\mathfrak{B}_n$ are equivalent, we need to show that there are natural isomorphisms $\overline{K}\overline{H} \longrightarrow\Id_{\overline{\mathfrak{N}}_n}$ and $\overline{H}\overline{K} \longrightarrow\Id_{\mathfrak{B}_n}$, where the $\Id_{\overline{\mathfrak{N}}_n}$ and $\Id_{\mathfrak{B}_n}$ are the respective identity functors on $\overline{\mathfrak{N}}_n$ and $\mathfrak{B}_n$. We first show that there is a natural isomorphism between functors $\overline{K}\overline{H}$ and $\Id_{\overline{\mathfrak{N}}_n}$. Since $\overline{h}: \mathcal{O}(\overline{\mathfrak{N}}_n) \longrightarrow\mathcal{O}(\mathfrak{B}_n)$ and $\overline{k}: \mathcal{O}(\mathfrak{B}_n)\longrightarrow\mathcal{O}(\overline{\mathfrak{N}}_n) $ are isomorphisms according to Lemma \ref{lemma:iso}, the composition of induced functions $\overline{k}\circ\overline{h}$ is the identity on $\mathcal{O}(\overline{\mathfrak{N}}_n)$. Therefore, for each $\mor\left(\overline{N}^\mathcal{X}_x,\overline{N}^\mathcal{X}_y\right)\in\mathcal{M}(\mathfrak{N}_n)$, we have a morphism $\beta_x: \overline{K}\overline{H}\left(\overline{N}^\mathcal{X}_x\right) \longrightarrow \Id_{\overline{\mathfrak{N}}_n}\left(\overline{N}^\mathcal{X}_x\right)$, $\beta_y: \overline{K}\overline{H}\left(\overline{N}^\mathcal{X}_y\right) \longrightarrow \Id_{\overline{\mathfrak{N}}_n}\left(\overline{N}^\mathcal{X}_y\right)$, $\phi^\prime \in\mathcal{M}(\overline{\mathfrak{N}}_n)$ sending the object $\overline{K}\overline{H} \left(\overline{N}^\mathcal{X}_x\right)$ to the object $\overline{K}\overline{H} \left(\overline{N}^\mathcal{X}_y\right)$, and $\psi^\prime \in\mathcal{M}(\mathfrak{B}_n)$ sending the object $\Id_{\overline{\mathfrak{N}}_n}\left(\overline{N}^\mathcal{X}_x\right) $ to the object $\Id_{\overline{\mathfrak{N}}_n}\left(\overline{N}^\mathcal{X}_y\right) $, such that
\begin{align*}
\beta_y \circ \phi^\prime \left(\overline{K}\overline{H} \left(\overline{N}^\mathcal{X}_x\right)\right)
 & = \beta_y \circ\overline{K}\overline{H} \left(\overline{N}^\mathcal{X}_y\right)\\
&=\Id_{\overline{\mathfrak{N}}_n}\left(\overline{N}^\mathcal{X}_y\right),\\
\psi^\prime\circ \beta_x \left( \overline{K}\overline{H} \left(\overline{N}^\mathcal{X}_x\right)\right)
&=\psi^\prime\circ  \Id_{\overline{\mathfrak{N}}_n}\left(\overline{N}^\mathcal{X}_x\right)\\
&=\Id_{\overline{\mathfrak{N}}_n}\left(\overline{N}^\mathcal{X}_y\right).
\end{align*}
Therefore, the diagram 
\[
\begin{tikzcd}[column sep=huge,row sep=huge]
\overline{K}\overline{H}\left(\overline{N}^\mathcal{X}_x\right) \arrow[r, "\phi^\prime"] \arrow[d,swap,"\beta_x"] &\overline{K}\overline{H}\left(\overline{N}^\mathcal{X}_y\right) \arrow[d,shift left=.75ex,"\beta_y"] \\
\Id_{\overline{\mathfrak{N}}_n}\left(\overline{N}^\mathcal{X}_x\right) \arrow[r,"\psi^\prime"] & \Id_{\overline{\mathfrak{N}}_n}\left(\overline{N}^\mathcal{X}_y\right).
\end{tikzcd}
\]
commutes.

Similarly, there are natural isomorphisms 
$\overline{H}\overline{K} \longrightarrow\Id_{\mathfrak{B}_n}$. For each $\mor\left(\mathcal{B}^X_x,\mathcal{B}^X_y\right)\in\mathcal{M}(\mathfrak{B}_n)$, by Lemma \ref{lemma:iso}, the composition of the induced functions $\overline{h}\circ\overline{k}$ is the identity on $\mathcal{O}(\mathfrak{B}_n)$. Therefore, we have morphisms $\gamma_x: \overline{H}\overline{K}\left(\mathcal{B}^X_x\right) \longrightarrow \Id_{\mathfrak{B}_n}\left(\mathcal{B}^X_x\right)$
and
 $\gamma_y: \overline{H}\overline{K}\left(\mathcal{B}^X_y\right) \longrightarrow \Id_{\mathfrak{B}_n}\left(\mathcal{B}^X_y\right)$, such that
\begin{align*}
\gamma_y \circ \phi^{\prime\prime} \left(\overline{H}\overline{K} \left(\mathcal{B}^X_x\right)\right)
 & = \gamma_y \circ\overline{H}\overline{K} \left(\mathcal{B}^X_y\right)\\
&=\Id_{\mathfrak{B}_n}\left(\mathcal{B}^X_y\right),\\
\psi^{\prime\prime}\circ \gamma_x \left( \overline{H}\overline{K} \left(\mathcal{B}^X_x\right)\right)
&=\psi^{\prime\prime}\circ  \Id_{\mathfrak{B}_n}\left(\mathcal{B}^X_x\right)\\
&=\Id_{\mathfrak{B}_n}\left(\mathcal{B}^X_y\right),
\end{align*}
  for a morphism $\phi^{\prime\prime} \in\mathcal{M}(\mathfrak{B}_n)$ sending the object  $\overline{H}\overline{K}(\mathcal{B}^X_x)$ to the object $\overline{H}\overline{K}(\mathcal{B}^X_y)$
 and  a morphism 
 $\psi^{\prime
\prime} \in\mathcal{M}(\overline{\mathfrak{N}}_n)$ sending the object $\Id_{\mathfrak{B}_n}(\mathcal{B}^X_x) $ to the object $\Id_{\mathfrak{B}_n}(\mathcal{B}^X_y) $.
Therefore, the diagram 
\[
\begin{tikzcd}[column sep=huge,row sep=huge]
\overline{H}\overline{K}(\mathcal{B}^X_x) \arrow[r, "\phi^{\prime\prime}"] \arrow[d,swap,"\gamma_x"] &\overline{H}\overline{K}(\mathcal{B}^X_y) \arrow[d,shift left=.75ex,"\gamma_y"] \\
\Id_{\mathfrak{B}_n}(\mathcal{B}^X_x) \arrow[r,"\psi^{\prime\prime}"] & \Id_{\mathfrak{B}_n}(\mathcal{B}^X_y).
\end{tikzcd}
\]
commutes.
Thus, we conclude that the category $\overline{\mathfrak{N}}_n$ and the category $\mathfrak{B}_n$ are equivalent.
\end{proof}

\section{Conclusion}

From the ever-changing genetic technology, new understandings of gene expression and genetic engineering emerge from day to day; however, the models of genetic evolution remain to be two classical models: the model of phylogenetic trees proposed by Darwin in 1859 to represent the evolutionary relationships of reproducing individuals \cite{darwin1859origin} and the model of phylogenetic networks proposed by Grant in 1971 to represent non-tree-like relationships such as species hybridization, bacterial gene transfer, and homologous recombination \cite{grant1971}. 
Therefore, it is of paramount importance to provide a rigorous definition of these spaces and investigate their connection with classical branches of mathematics.
 
In this article, we focused on the graph theoretical aspects of phylogenetic trees and networks and their connection to stable curves in  $\overline{\mathfrak{M}}_{0,n}(\mathbb{C})$, $\overline{\mathfrak{M}}_{g,n}(\mathbb{C})$, and  $\overline{\mathfrak{M}}_{0,n}(\mathbb{R})$. We introduced building blocks of evolutionary moduli spaces and the dual intersection complex of the moduli spaces of stable curves, and the equivalence of categories between these spaces. We defined the category of trees $\ft_n$, the category of stable curves with $n$ marked points $\fc_n$, and the category of tree topologies $\overline{\mathfrak{T}}_n$. We proved that there is a functor from the category $\fc_n$ to the category of tree topologies $\overline{\mathfrak{T}}_n$.
We defined the category of the partition space of genus 0 stable curves with $n$ marked points $\ofc_n$ and showed that there is a functor from $\mathfrak{T}_n$
to $\ofc_n$. Then we proved that the category $\overline{\fc}_n$ is equivalent to the category $\overline{\mathfrak{T}}_n$.  Another important result we reveal between the space of phylogenetic trees and stable curves is that the space of stable curves and tree space are the dual of each other such that a tree with $k$ internal edges and $n$ nodes corresponds to the moduli space with $n-3-k$ complex parameters in $\overline{\mathfrak{M}}_{0,n}(\mathbb{C})$. We also showed that the dual intersection complex of boundary divisors of $\overline{\mathfrak{M}}_{0,n}(\mathbb{C})$ as an abstract simplicial complex is isomorphic to $\pbhv_{n-1}$ as an abstract simplicial complex. Consequently, their respective geometric realizations are homeomorphic.

For the space of phylogenetic networks, we provided a formal definition of phylogenetic networks and the space in which phylogenetic networks reside. We first investigated the dual graphs of phylogenetic networks that are associated with stable curves for arbitrary genus, and we proved the seemingly distant relation that phylogenetic networks in $\csn_n$ always correspond to a stable curve in $\overline{\mathfrak{M}}_{g,n}$. Then we defined the category of network $\mathfrak{N}_n$ and the category of boundary divisors of $\overline{\mathfrak{M}}_{0,n}(\mathbb{R})$, denoted as $\mathfrak{B}_n$. We showed that there is a functor which maps from the category of networks $\mathfrak{N}_n$ onto the category of boundary divisors $\mathfrak{B}_n$. The category of networks and the category of boundary divisors are connected through network topologies. We defined the category of network topologies $\overline{\mathfrak{N}}_n$ and showed that the category $\overline{\mathfrak{N}}_n$ and the category $\mathfrak{B}_n$ are equivalent.

Regardless of the results in this paper, there is still much to be understood.  For example, the space of phylogenetic trees and the space of phylogenetic networks appear in the work of Hacking, Keel, and Tevelev but based there on the real algebraic geometry of del Pezzo surfaces and root systems \cite{hacking2009stable}.  It is thought that varying the root system will yield a spectrum of spaces analogous to $\bhv_n$ and $\csn_n$, and that this connection with root systems should shed light on the geometric properties of phylogenetic spaces \cite{devadoss2017space}.  By way of a second example: It may well be worthwhile to explore the combinatorial properties of phylogenetic spaces in more detail, for there are still open questions in this direction.

\bibliographystyle{plain}
\bibliography{bib} 

\begin{thebibliography}{10}

\bibitem{arbarello2011geometry}
Enrico Arbarello, Maurizio Cornalba, and Phillip Griffiths.
\newblock {\em Geometry of algebraic curves: volume II with a contribution by
  Joseph Daniel Harris}, volume 268.
\newblock Springer Science \& Business Media, 2011.

\bibitem{ardila2012geodesics}
Federico Ardila, Megan Owen, and Seth Sullivant.
\newblock Geodesics in cat (0) cubical complexes.
\newblock {\em Advances in Applied Mathematics}, 48(1):142--163, 2012.

\bibitem{baez2017operads}
John~C Baez and Nina Otter.
\newblock Operads and phylogenetic trees.
\newblock {\em Theory and Applications of Categories}, 32(40):1397--1453, 2017.

\bibitem{billera2001geometry}
Louis~J Billera, Susan~P Holmes, and Karen Vogtmann.
\newblock Geometry of the space of phylogenetic trees.
\newblock {\em Advances in Applied Mathematics}, 27(4):733--767, 2001.

\bibitem{bryant2002neighbornet}
David Bryant and Vincent Moulton.
\newblock Neighbornet: An agglomerative method for the construction of planar
  phylogenetic networks.
\newblock In {\em WABI}, volume~2, pages 375--391. Springer, 2002.

\bibitem{darwin1859origin}
Charles Darwin.
\newblock {\em On the origin of the species by natural selection}.
\newblock Murray, 1859.

\bibitem{deligne1969irreducibility}
Pierre Deligne and David Mumford.
\newblock The irreducibility of the space of curves of given genus.
\newblock {\em Publications Math{\'e}matiques de l'IHES}, 36:75--109, 1969.

\bibitem{deligne1972irreducibility}
Pierre Deligne and David Mumford.
\newblock The irreducibility of the space of curves of given genus.
\newblock {\em Matematika}, 16(3):13--53, 1972.

\bibitem{devadoss2010diagonalizing}
Satyan~L Devadoss and Jack Morava.
\newblock Diagonalizing the genome i: navigation in tree spaces.
\newblock {\em arXiv preprint arXiv:1009.3224}, 2010.

\bibitem{devadoss2017space}
Satyan~L Devadoss and Samantha Petti.
\newblock A space of phylogenetic networks.
\newblock {\em SIAM Journal on Applied Algebra and Geometry}, 1(1):683--705,
  2017.

\bibitem{doolittle1998you}
W~Ford Doolittle.
\newblock You are what you eat: a gene transfer ratchet could account for
  bacterial genes in eukaryotic nuclear genomes.
\newblock {\em Trends in Genetics}, 14(8):307--311, 1998.

\bibitem{Edelman594}
Nathaniel~B. Edelman, Paul~B. Frandsen, Michael Miyagi, Bernardo Clavijo, John
  Davey, Rebecca~B. Dikow, Gonzalo Garc{\'\i}a-Accinelli, Steven~M.
  Van~Belleghem, Nick Patterson, Daniel~E. Neafsey, Richard Challis, Sujai
  Kumar, Gilson R.~P. Moreira, Camilo Salazar, Mathieu Chouteau, Brian~A.
  Counterman, Riccardo Papa, Mark Blaxter, Robert~D. Reed, Kanchon~K.
  Dasmahapatra, Marcus Kronforst, Mathieu Joron, Chris~D. Jiggins, W.~Owen
  McMillan, Federica Di~Palma, Andrew~J. Blumberg, John Wakeley, David Jaffe,
  and James Mallet.
\newblock Genomic architecture and introgression shape a butterfly radiation.
\newblock {\em Science}, 366(6465):594--599, 2019.

\bibitem{fulton1996notes}
W.~Fulton and R.~Pandharipande.
\newblock Notes on stable maps and quantum cohomology, 1996.

\bibitem{grant1971}
Verne Grant.
\newblock {\em {P}lant {S}peciation}.
\newblock Columbia University Press, New York, 2nd edition, 1981.

\bibitem{hacking2009stable}
Paul Hacking, Se{\'a}n Keel, and Jenia Tevelev.
\newblock Stable pair, tropical, and log canonical compactifications of moduli
  spaces of del pezzo surfaces.
\newblock {\em Inventiones mathematicae}, 178(1):173--227, 2009.

\bibitem{harris2006moduli}
Joe Harris and Ian Morrison.
\newblock {\em Moduli of curves}, volume 187.
\newblock Springer Science \& Business Media, 2006.

\bibitem{hartshorne2013algebraic}
Robin Hartshorne.
\newblock {\em Algebraic geometry}, volume~52.
\newblock Springer Science \& Business Media, 2013.

\bibitem{allen2002cambridge}
Allen Hatcher.
\newblock {\em Algebraic Topology}.
\newblock Cambridge University Press, Cambridge, 2002.

\bibitem{hotopp2011horizontal}
Julie C~Dunning Hotopp.
\newblock Horizontal gene transfer between bacteria and animals.
\newblock {\em Trends in genetics}, 27(4):157--163, 2011.

\bibitem{huson2006application}
Daniel~H Huson and David Bryant.
\newblock Application of phylogenetic networks in evolutionary studies.
\newblock {\em Molecular biology and evolution}, 23(2):254--267, 2006.

\bibitem{huson2010phylogenetic}
Daniel~H Huson, Regula Rupp, and Celine Scornavacca.
\newblock {\em Phylogenetic networks: concepts, algorithms and applications}.
\newblock Cambridge University Press, 2010.

\bibitem{jain2003horizontal}
Ravi Jain, Maria~C Rivera, Jonathan~E Moore, and James~A Lake.
\newblock Horizontal gene transfer accelerates genome innovation and evolution.
\newblock {\em Molecular biology and evolution}, 20(10):1598--1602, 2003.

\bibitem{keel1992intersection}
Se{\'a}n Keel.
\newblock Intersection theory of moduli space of stable n-pointed curves of
  genus zero.
\newblock {\em Transactions of the American Mathematical Society}, pages
  545--574, 1992.

\bibitem{kozak2018genome}
Krzysztof~M Kozak, Owen McMillan, Mathieu Joron, and Christopher~D Jiggins.
\newblock Genome-wide admixture is common across the heliconius radiation.
\newblock {\em BioRxiv}, page 414201, 2018.

\bibitem{lamichhaney2015evolution}
Sangeet Lamichhaney, Jonas Berglund, Markus~S{\"a}llman Alm{\'e}n, Khurram
  Maqbool, Manfred Grabherr, Alvaro Martinez-Barrio, Marta Promerov{\'a},
  Carl-Johan Rubin, Chao Wang, Neda Zamani, et~al.
\newblock Evolution of darwin’s finches and their beaks revealed by genome
  sequencing.
\newblock {\em Nature}, 518(7539):371--375, 2015.

\bibitem{may1997definitions}
J~Peter May.
\newblock Definitions: operads, algebras and modules.
\newblock {\em Contemporary Mathematics}, 202:1--8, 1997.

\bibitem{may1999concise}
J~Peter May.
\newblock {\em A concise course in algebraic topology}.
\newblock University of Chicago press, 1999.

\bibitem{moerdijk2007dendroidal}
Ieke Moerdijk and Ittay Weiss.
\newblock Dendroidal sets.
\newblock {\em Algebraic \& Geometric Topology}, 7(3):1441--1470, 2007.

\bibitem{owen2011computing}
Megan Owen.
\newblock Computing geodesic distances in tree space.
\newblock {\em SIAM Journal on Discrete Mathematics}, 25(4):1506--1529, 2011.

\bibitem{owen2011fast}
Megan Owen and J~Scott Provan.
\newblock A fast algorithm for computing geodesic distances in tree space.
\newblock {\em IEEE/ACM Transactions on Computational Biology and
  Bioinformatics (TCBB)}, 8(1):2--13, 2011.

\bibitem{pease2016phylogenomics}
James~B Pease, David~C Haak, Matthew~W Hahn, and Leonie~C Moyle.
\newblock Phylogenomics reveals three sources of adaptive variation during a
  rapid radiation.
\newblock {\em PLoS Biology}, 14(2), 2016.

\bibitem{saitou1987neighbor}
Naruya Saitou and Masatoshi Nei.
\newblock The neighbor-joining method: a new method for reconstructing
  phylogenetic trees.
\newblock {\em Molecular biology and evolution}, 4(4):406--425, 1987.

\bibitem{speyer2004tropical}
David Speyer and Bernd Sturmfels.
\newblock The tropical grassmannian.
\newblock {\em Advances in Geometry}, 4(3):389--411, 2004.

\bibitem{wu2019comparison}
Yingying Wu.
\newblock {\em Comparison Theorems of Phylogenetic Spaces and Algebraic Fans}.
\newblock Dissertation, The University of Texas at Austin, 2019.

\bibitem{zairis2014moduli}
Sakellarios Zairis, Hossein Khiabanian, Andrew~J Blumberg, and Raul Rabadan.
\newblock Moduli spaces of phylogenetic trees describing tumor evolutionary
  patterns.
\newblock In {\em International Conference on Brain Informatics and Health},
  pages 528--539. Springer, 2014.

\bibitem{zairis2016genomic}
Sakellarios Zairis, Hossein Khiabanian, Andrew~J Blumberg, and Raul Rabadan.
\newblock Genomic data analysis in tree spaces.
\newblock {\em arXiv preprint arXiv:1607.07503}, 2016.

\end{thebibliography}
\index{Bibliography@\emph{Bibliography}}

\vskip 12mm

\vskip 12mm

\begin{tabular}{ ll }
 Y. Wu: & Center of Mathematical Sciences and Applications\\
 & Harvard University\\ 
& Cambridge, MA, 02138\\
&	\url{ywu@cmsa.fas.harvard.edu}\\
 &\\
S.-T. Yau: &Department of Mathematics\\
& Harvard University\\
& Cambridge, MA, 02138\\
 &\url{yau@math.harvard.edu}\\
\end{tabular}
\\

\end{document}